\documentclass[12pt,leqno]{article}
\usepackage{amsmath, amsthm, amsfonts, amssymb, color, enumerate}
\usepackage{mathrsfs}
\newtheorem{thm}{Theorem}[section]
\newtheorem{cor}[thm]{Corollary}
\newtheorem{lem}[thm]{Lemma}
\newtheorem{prp}[thm]{Proposition}

\theoremstyle{definition}
\newtheorem{defn}{Definition}[section]
\newcommand{\scr}[1]{\mathscr #1}
\definecolor{wco}{rgb}{0.5,0.2,0.3}

\numberwithin{equation}{section} \theoremstyle{remark}

\newcommand{\ua}{\uparrow}

\title{
{\bf  Degenerate SDE with H\"older-Dini Drift and Non-Lipschitz Noise Coefficient}
\footnote{FW is supported in part by
NNSFC (11131003, 11431014), the 985 project and the Laboratory of Mathematical and  Complex Systems, XZ is supported partly by
NNSFC (11271294, 11325105).}}
\author{
{\bf Feng-Yu Wang$^{a),b)}$ and Xicheng Zhang$^{c)}$ }\\
\footnotesize{a) School of Mathematical Sciences,
Beijing Normal University, Beijing 100875, China}\\
 \footnotesize{b) Department of Mathematics,
Swansea University, Singleton Park, SA2 8PP, UK}\\
\footnotesize{c) School of Mathematics and Statistics,
Wuhan University, Wuhan 430072, China}}

\begin{document}
\numberwithin{equation}{section}
\def\theequation{\arabic{section}.\arabic{equation}}

\newcommand{\be}{\begin{eqnarray}}
\newcommand{\ee}{\end{eqnarray}}
\newcommand{\ce}{\begin{eqnarray*}}
\newcommand{\de}{\end{eqnarray*}}
\newtheorem{theorem}{Theorem}[section]
\newtheorem{lemma}[theorem]{Lemma}
\newtheorem{remark}[theorem]{Remark}
\newtheorem{definition}[theorem]{Definition}
\newtheorem{proposition}[theorem]{Proposition}
\newtheorem{Examples}[theorem]{Example}
\newtheorem{corollary}[theorem]{Corollary}

\def\Re{{\mathrm{Re}}}
\def\Im{{\mathrm{Im}}}
\def\var{{\mathrm{var}}}
\def\HS{{\mathrm{\tiny HS}}}
\def\eps{\varepsilon}
\def\t{\tau}
\def\e{\mathrm{e}}
\def\vr{\varrho}
\def\th{\theta}
\def\a{\alpha}
\def\om{\omega}
\def\Om{\Omega}
\def\v{\mathrm{v}}
\def\u{\mathbf{u}}
\def\w{\mathrm{w}}
\def\p{\partial}
\def\d{\delta}
\def\g{\gamma}
\def\l{\lambda}
\def\la{\langle}
\def\ra{\rangle}
\def\[{{\Big[}}
\def\]{{\Big]}}
\def\<{{\langle}}
\def\>{{\rangle}}
\def\({{\Big(}}
\def\){{\Big)}}
\def\bz{{\mathbf{z}}}
\def\by{{\mathbf{y}}}
\def\bx{{\mathbf{x}}}
\def\tr{\mathrm {tr}}
\def\W{{\mathcal W}}
\def\Ric{{\rm Ricci}}
\def\osc{{\rm osc}}
\def\Cap{\mbox{\rm Cap}}
\def\sgn{\mbox{\rm sgn}}
\def\mathcalV{{\mathcal V}}
\def\Law{{\mathord{{\rm Law}}}}
\def\dif{{\mathord{{\rm d}}}}
\def\dis{{\mathord{{\rm \bf d}}}}
\def\Hess{{\mathord{{\rm Hess}}}}
\def\min{{\mathord{{\rm min}}}}
\def\Vol{\mathord{{\rm Vol}}}
\def\bbbn{{\rm I\!N}}
\def\no{\nonumber}
\def\={&\!\!=\!\!&}
\def\bt{\begin{thm}}
\def\et{\end{thm}}
\def\bl{\begin{lem}}
\def\el{\end{lem}}
\def\br{\begin{remark}}
\def\er{\end{remark}}
\def\bd{\begin{definition}}
\def\ed{\end{definition}}
\def\bp{\begin{prp}}
\def\ep{\end{prp}}
\def\bc{\begin{cor}}
\def\ec{\end{cor}}
\def\bx{\begin{Examples}}
\def\ex{\end{Examples}}

\def\cA{{\mathcal A}}
\def\cB{{\mathcal B}}
\def\cC{{\mathcal C}}
\def\cD{{\mathcal D}}
\def\cE{{\mathcal E}}
\def\cF{{\mathcal F}}
\def\cG{{\mathcal G}}
\def\cH{{\mathcal H}}
\def\cI{{\mathcal I}}
\def\cJ{{\mathcal J}}
\def\cK{{\mathcal K}}
\def\cL{{\mathcal L}}
\def\cM{{\mathcal M}}
\def\cN{{\mathcal N}}
\def\cO{{\mathcal O}}
\def\cP{{\mathcal P}}
\def\cQ{{\mathcal Q}}
\def\cR{{\mathcal R}}
\def\cS{{\mathcal S}}
\def\cT{{\mathcal T}}
\def\cU{{\mathcal U}}
\def\cV{{\mathcal V}}
\def\cW{{\mathcal W}}
\def\cX{{\mathcal X}}
\def\cY{{\mathcal Y}}
\def\cZ{{\mathcal Z}}

\def\mA{{\mathbb A}}
\def\mB{{\mathbb B}}
\def\mC{{\mathbb C}}
\def\mD{{\mathbb D}}
\def\mE{{\mathbb E}}
\def\mF{{\mathbb F}}
\def\mG{{\mathbb G}}
\def\mH{{\mathbb H}}
\def\mI{{\mathbb I}}
\def\mJ{{\mathbb J}}
\def\mK{{\mathbb K}}
\def\mL{{\mathbb L}}
\def\mM{{\mathbb M}}
\def\mN{{\mathbb N}}
\def\mO{{\mathbb O}}
\def\mP{{\mathbb P}}
\def\mQ{{\mathbb Q}}
\def\mR{{\mathbb R}}
\def\mS{{\mathbb S}}
\def\mT{{\mathbb T}}
\def\mU{{\mathbb U}}
\def\mV{{\mathbb V}}
\def\mW{{\mathbb W}}
\def\mX{{\mathbb X}}
\def\mY{{\mathbb Y}}
\def\mZ{{\mathbb Z}}

\def\bB{{\mathbf B}}
\def\bP{{\mathbf P}}
\def\bQ{{\mathbf Q}}
\def\bE{{\mathbf E}}
\def\1{{\mathbf{1}}}

\def\sA{{\mathscr A}}
\def\sB{{\mathscr B}}
\def\sC{{\mathscr C}}
\def\sD{{\mathscr D}}
\def\sE{{\mathscr E}}
\def\sF{{\mathscr F}}
\def\sG{{\mathscr G}}
\def\sH{{\mathscr H}}
\def\sI{{\mathscr I}}
\def\sJ{{\mathscr J}}
\def\sK{{\mathscr K}}
\def\sL{{\mathscr L}}
\def\sM{{\mathscr M}}
\def\sN{{\mathscr N}}
\def\sO{{\mathscr O}}
\def\sP{{\mathscr P}}
\def\sQ{{\mathscr Q}}
\def\sR{{\mathscr R}}
\def\sS{{\mathscr S}}
\def\sT{{\mathscr T}}
\def\sU{{\mathscr U}}
\def\sV{{\mathscr V}}
\def\sW{{\mathscr W}}
\def\sX{{\mathscr X}}
\def\sY{{\mathscr Y}}
\def\sZ{{\mathscr Z}}\def\E{\mathbb E}

\def\fM{{\mathfrak M}}
\def\fA{{\mathfrak A}}

\def\geq{\geqslant}
\def\leq{\leqslant}
\def\ge{\geqslant}
\def\le{\leqslant}

\def\c{\mathord{{\bf c}}}
\def\div{\mathord{{\rm div}}}
\def\iint{\int\!\!\!\int}

\def\sb{{\mathfrak b}}

\def\bH{{\mathbf H}}
\def\bW{{\mathbf W}}
\def\bP{{\mathbf P}}
\def\bA{{\mathbf A}}\def\si{\sigma}
\def\bT{{\mathbf T}} \def\R{\mathbb R}\def\nn{\nabla}
\def\ff{\frac} \def\R{\mathbb R}  \def\ff{\frac} \def\ss{\sqrt} \def\B{\scr B}
\def\N{\mathbb N} \def\kk{\kappa} \def\m{{\bf m}}
\def\dd{\delta} \def\DD{\Dd} \def\vv{\varepsilon} \def\rr{\rho}
\def\<{\langle} \def\>{\rangle} \def\GG{\Gamma} \def\gg{\gamma}
  \def\nn{\nabla} \def\pp{\partial} \def\EE{\scr E}
\def\d{\text{\rm{d}}} \def\bb{\beta} \def\aa{\alpha} \def\D{\scr D}
  \def\si{\sigma} \def\ess{\text{\rm{ess}}}
\def\beg{\begin} \def\beq{\begin{equation}}  \def\F{\scr F}
\def\Ric{\text{\rm{Ric}}} \def\Hess{\text{\rm{Hess}}}
\def\e{\text{\rm{e}}} \def\ua{\underline a} \def\OO{\Omega}  \def\oo{\omega}
 \def\tt{\tilde} \def\Ric{\text{\rm{Ric}}}
\def\cut{\text{\rm{cut}}} \def\P{\mathbb P} \def\ifn{I_n(f^{\bigotimes n})}
\def\C{\scr C}      \def\aaa{\mathbf{r}}     \def\r{r}
\def\gap{\text{\rm{gap}}} \def\prr{\pi_{{\bf m},\varrho}}  \def\r{\mathbf r}
\def\Z{\mathbb Z} \def\vrr{\varrho} \def\l{\lambda}\def\ppp{\preceq}
\def\L{\scr L}\def\Tt{\tt} \def\TT{\tt}\def\II{\mathbb I}\def\ll{\lambda} \def\LL{\Lambda}
\allowdisplaybreaks
\maketitle
  \begin{abstract} The existence-uniqueness and stability of strong  solutions are proved for   a class of degenerate stochastic
differential equations,  where the noise coeffcicient might be non-Lipschitz, and the drift is locally Dini continuous in  the component
with noise (i.e. the second component) and locally H\"older-Dini continuous of order   $\ff 2 3$ in
the first component. Moreover, the weak uniqueness is proved under  weaker conditions on the noise coefficient.
Furthermore, if the noise coefficient is  $C^{1+\vv}$ for some $\vv>0$ and the drift is   H\"older continuous  of order $\aa\in (\ff 2 3,1)$ in the first component and order $\bb\in(0,1) $ in the second,  the solution forms a $C^1$-stochastic diffeormorphism flow.
To prove these results, we present some new characterizations of H\"older-Dini space by using the heat semigroup and slowly varying functions.

\end{abstract} \noindent


\rm

\section{Introduction}

Consider the following ordinary differential equation (abbreviated as ODE):
$$
\dot x(t)=b(x(t)),\ x(0)=x_0.
$$
It is   classical   that the equation is well-posed for Lipschitz $b$ but usually ill-posed if $b$ is only H\"older continuous.
For instance, for $b(x):=|x|^{\alpha}$ with $\alpha\in(0,1)$ and $x_0=0$,
 the above ODE has two solutions: $x(t)\equiv0$ and $x(t)=(1-\alpha)t^{1/(1-\alpha)}, t\geq 0$.
However,  if the above ODE is perturbed by a strong enough   noise  (e.g. the Browian motion),
the equation might be well-posed for very singular $b$. For instance,   consider the following SDE on $\mR^d$:
$$
\dif X_t=b_t(X_t)\dif t+\sigma \dif W_t,\ \ X_0=x,
$$
where $W_t$ is a $d$-dimensional standard Brownian motion on some probability space $(\Omega,\sF,\mP)$, $\sigma$ is an invertible matrix. If $b$ is a bounded measurable function,
Veretennikov \cite{Ve} proved that the above SDE admits a unique strong solution, which extended an earlier result of Zvonkin \cite{Zv} in the case of $d=1$.
More recent results about the above SDE can be found in \cite{Fe-Fl, Kr-Ro, Zh0} and references therein for further development in this direction.

It is worthy noticing that all the well-posedness results mentioned above are done only for the {\it time-white} noise, which means that the noise is a distribution of the time variable.
In this work, we are concerning with the following problem: Is it possible to prove the well-posedness of the ODE with singular $b$
perturbed by  an absolutely continuous Gaussian process?    More concretely,   consider the following random ODE:
\begin{align}
\dif X_t=[b_t(X_t)+\sigma  W_t]\dif t,\ \ X_0=x.\label{SDE0}
\end{align}
We aim to find minimal conditions on $b$ and $\si$ ensuring  the well-posedness of this random ODE. By regarding $X_t$ as the first component process $X_t^{(1)}$ and introducing $X_t^{(2)}:= \si W_t$, this problem is reduced to
 the study of   the following more general   degenerate SDE for $X_t:=(X_t^{(1)}, X_t^{(2)})$ on
$\R^{d_1+d_2}=\R^{d_1}\times\R^{d_2}$:
\begin{align}
\dif X_t=b_t(X_t)\dif t+(0,\sigma_t(X_t)\dif W_t), \ \ X_0=x=(x^{(1)}, x^{(2)})\in \mR^{d_1+d_2},\label{SDE}
\end{align}
where, for $\R_+:=(0,\infty),$ the maps $\si: \R_+\times\R^{d_1+d_2}\to \R^{d_2}\otimes\R^{d_2}$ and $ b=(b^{(1)}, b^{(2)}): \R_+\times \R^{d_1+d_2} \to \R^{d_1+d_2}$
are measurable and locally bounded. This model is known as the stochastic Hamiltonian system with potential $H$ if $b=\nn H$, which includes the kinetic Fokker-Planck equation as a typical example (see \cite{V}).

In the following, we will use $\nn^{(1)}$ and $\nn^{(2)}$ to denote the gradient operators on the first space $\R^{d_1}$
and the second space $\R^{d_2}$ respectively. Thus, for every $(t,x)\in \R_+\times \R^{d_1+d_2}$,
$\nn^{(2)}b_t^{(1)}(x)\in \R^{d_2}\otimes \R^{d_1}$ with $(\nn^{(2)}b_t^{(1)}(x))h:= \nn^{(2)}_h b_t^{(1)}(x)\in \R^{d_1}, h\in \R^{d_2}.$
By It\^o's formula,  the infinitesimal generator associated to \eqref{SDE} is given by
\begin{align}
\sL^{\Sigma,b}_tu=\mathrm{tr}\big(\Sigma_t\cdot\nabla^{(2)}\nabla^{(2)}u\big)+b_t\cdot\nabla u,\label{LG}
\end{align}
where $\Sigma_t(x):= \tfrac{1}{2}\si_t(x)\si_t^*(x)$ and tr$(\cdot)$ denotes the trace of a matrix.

Let $|\cdot|$ denote the Euclidiean norm and let $\|\cdot\|$ denote the operator norm. We introduce below the notion of H\"older-Dini continuity.

\begin{defn}
An increasing function $\phi:\R_+\to\R_+$ is called a {\it Dini function}  if
 \begin{align}
 \int_0^1 \ff{\phi(t)}t \d t<\infty.\label{YY3}
 \end{align}
A measurable function $\phi:\R_+\to\R_+$  is called a {\it slowly varying function} at zero if for any $\lambda>0$,
\begin{align}
\lim_{t\to 0}\frac{\phi(\lambda t)}{t}=1.\label{YY2}
\end{align}A function $f$ on the Euclidiean space is called {\it H\"older-Dini continuous} of order $\aa\in [0,1)$ if
$$|f(x)-f(y)|\le |x-y|^\aa\phi(|x-y|),\ \ |x-y|\le 1$$ holds for some Dini function $\phi$,
and is called {\it Dini-continuous} if this condition holds for $\aa=0$.
 \end{defn}

Let $\sD_0$ be the set of all Dini functions, and $\sS_0$ the set of all slowly varying functions that are bounded from $0$ and $\infty$ on $[\eps,\infty)$ for any $\eps>0$.
Notice that the typical examples in $\sD_0\cap\sS_0$ are
$\phi(t):= (\log (1+t^{-1}))^{-\beta}$ for $\beta>1$.

Roughly speaking, for the existence and uniqueness of the solutions to \eqref{SDE}, we will need   $b^{(1)}(\cdot, x^{(2)})$
and $b^{(2)}(\cdot, x^{(2)})$ and $\nn^{(2)} b^{(2)}(x^{(1)},\cdot)$ with fixed $x^{(2)}$ to be locally H\"older-Dini continuous of order   $\ff 2 3$,  and $b^{(2)}(x^{(1)},\cdot)$ with fixed $x^{(1)}$ to be merely Dini continuous. These coincide  with the continuity conditions used in \cite{Wa-Zh} for infinite-dimensional degenerate systems with linear
 $b^{(1)}$.

 Moreover, it is known that \eqref{SDE} is well-posed if   $\si$ and $b$ are $``$almostly Lipschitz continuous'', see e.g. \cite{Ya-Wa, FZ, SWY}. In this paper we show that,  under the above mentioned much weaker conditions on $b$,  such a non-Lipschitz condition on $\si$ still implies the well-posedness.  To  characterize this condition,    we introduce  the class
$$\C:= \bigg\{\gg\in C^1(\R_+;\R_+):\ \int_0^1\ff{1}{t\gg(t)}\d t=\infty,\ \liminf_{t\downarrow 0} \Big(\ff {\gg(t)} 4 +t\gg'(t)\Big)>0  \bigg\},$$
where   $\int_0^1\ff{1}{t\gg(t)}\d t=\infty$ is the key condition, and $\liminf_{t\downarrow 0} \big(\ff {\gg(t)} 4 +t\gg'(t)\big)>0$  comes from our calculations in  the present framework, which is weaker than the following condition used in \cite[Theorem B]{FZ}:
$$\lim_{t\downarrow 0} \gg(t)=\infty,\ \ \lim_{t\downarrow 0} \ff{t\gg'(t)}{\gg(t)}=0.$$
Typical functions in  $\C$  include $$\gg_1(t):= \log (1+t^{-1}), \ \gg_2(t):=\gg_1(t) \log\log (\e+t^{-1}),\
\gg_3(t):=\gg_2(t)  \log\log\log (\e^2+t^{-1})...$$

In the following four subsections, we state our main results on the weak solutions, the strong solutions, the stability of solutions with respect to coefficients, and the $C^1$-stochastic diffeormorphism flows respectively.

\subsection{Weak solutions}

 We introduce the following assumptions for some $\phi\in \D_0\cap\sS_0$ and some increasing function $C:\R_+\to\R_+$:

\beg{enumerate}
\item[{\bf (C1)}]   {\bf (Hypoellipticity)} $\si_t(x)$ and $[\nn^{(2)} b_t^{(1)}(x)][\nn^{(2)} b_t^{(1)}(x)]^*$ are invertible  with
\beg{align*}
  \| \nn^{(2)} b_t^{(1)}  \|_\infty +\big\|\big([\nn^{(2)} b_t^{(1)} ][\nn^{(2)} b_t^{(1)} ]^*\big)^{-1}\big\|_\infty+\|\si_t\|_\infty+ \|\si_t^{-1}\|_\infty \le C (t),\ \  t\ge 0.
\end{align*}
 \item[{\bf(C2)}] {\bf (Regularity of $b^{(1)}$)} For any $x,y\in \R^{d_1+d_2}$ with $|x-y|\le 1$ and $t\ge 0$,
\beg{equation*}\beg{split}
 &|b^{(1)}_t(x)-b^{(1)}_t(y)|\leq  C(t) |x^{(1)}-y^{(1)}|^{\ff 2 3}\phi(|x^{(1)}-y^{(1)}|),\ \quad \text{if} \ x^{(2)}=y^{(2)},\\
&\|\nabla^{(2)}b^{(1)}_t(x)-\nabla^{(2)}b^{(1)}_t(y)\|\leq C(t)\phi(|x^{(2)}-y^{(2)}|),\   \quad \text{if} \ x^{(1)}=y^{(1)}.
\end{split}\end{equation*}
\item[{\bf(C3)}] {\bf (Regularity of $b^{(2)},\si$)}
Either
\beq\label{LB1}\beg{cases}\beg{split}
 &|b^{(2)}_t(x)-b^{(2)}_t(y)| \leq C (t)\big\{|x^{(1)}-y^{(1)}|^{\ff 2 3}\phi (|x^{(1)}-y^{(1)}|) +\phi^{\ff 72} (|x^{(2)}-y^{(2)}|)\big\},\\
 &\|\si_t(x)-\si_t(y)\|\le C(t)|x-y|^{\ff 2 3} \phi(|x-y|),\ \ t\ge 0, |x-y|\le 1; \end{split}\end{cases}\end{equation} or for  $ t\ge 0, |x-y|\le 1,$ there hold $\|\nn^{(2)}\si_t\|_\infty\le C(t)$ and
\beq\label{LB2}\beg{cases}\beg{split}
&|b^{(2)}_t(x)-b^{(2)}_t(y)| \leq C (t)\big\{|x^{(1)}-y^{(1)}|^{\ff 2 3}\phi (|x^{(1)}-y^{(1)}|) +\phi  (|x^{(2)}-y^{(2)}|)\big\},\\
& \|\nn^{(2)} \si_t(x^{(1)},x^{(2)})-\nn^{(2)}\si_t(y^{(1)}, x^{(2)})\|\le C(t) |x^{(1)}-y^{(1)}|^{\ff 19} \phi(|x^{(1)}-y^{(1)}|),\\
&\|\si_t(x^{(1)},x^{(2)})-\si_t(y^{(1)}, x^{(2)})\|\le C(t)|x^{(1)}-y^{(1)}|^{\ff 2 3} \phi(|x^{(1)}-y^{(1)}|).\end{split}\end{cases}\end{equation}
  \end{enumerate}

Intuitively,   there should be a balance between the regularities of $b^{(2)}$ and   $\si$; that is, with a stronger condition on $\si$ we will only need a weaker regularity of $b^{(2)}$. Conditions \eqref{LB1} and \eqref{LB2}, as well as \eqref{LB1'} and \eqref{LB2'} below, are introduced in this spirit.

\beg{thm}\label{TW} Assume that
{\bf (C1)}--{\bf (C3)}  hold  for some $\phi\in \D_0\cap\sS_0$ and   increasing function $C: \R_+\to\R_+$. Then $\eqref{SDE}$ has a unique weak solution. \end{thm}

\paragraph{Remark 1.1.} In \cite{Me}, Menozzi showed that the weak uniqueness holds for \eqref{SDE}
under the assumptions that $\sigma$ is H\"older continuous
and $b$ is Lipschitz continuous.
In \cite{Pr}, Priola showed that there is a unique weak solution to \eqref{SDE} when $\sigma_t(x)=\sigma(x)$
is bounded continuous, $b^{(1)}(x)=x^{(2)}$ and $b^{(2)}(x)$ is bounded measurable. Although our assumptions on $b^{(2)}$ and $\sigma$
are stronger,  we allow   $b^{(1)}(x)$   to be merely H\"older-Dini continuous in $x^{(1)}$. In fact, this is the main source of the difficulty in our study, since due to the singularity of $b^{(1)}(x)$ in $x^{(1)}$ we have to carefully estimate the regularization of the noise transported from the second component to the first, see
  Lemma 3.1 below.

\subsection{Strong solutions}

By a localization argument, we will take the following local conditions on $\si$ and $b$.

\beg{enumerate} \item[{\bf (A)}\ \ ] For any $n\in\N$,  there exist a constant $C_n\in\R_+$,   some  $\phi_n\in \D_0\cap\sS_0$ and $\gg_n\in \C$
such that the following conditions hold for all $t\in [0,n]$:
\item[{\bf (A1)}]   {\bf (Hypoellipticity)} $\si_t(x)$ and $[\nn^{(2)} b_t^{(1)}(x)][\nn^{(2)} b_t^{(1)}(x)]^*$ are invertible and locally bounded with
\beg{align*}
  \sup_{x\in \R^{d_1+d^2}, |x^{(1)}|\le n}\big\| \big([\nn^{(2)} b_t^{(1)} ][\nn^{(2)} b_t^{(1)} ]^*\big)^{-1}\big\|(x) +\sup_{|x|\le n} \|\si_t^{-1}\|(x) \le C_n.
\end{align*}
 \item[{\bf(A2)}] {\bf (Regularity of $b^{(1)}$)} For any $x,y\in \R^{d_1+d_2}$ with $ |x|\lor |y|\le n$,
\beg{equation*}\beg{split}
 &|b^{(1)}_t(x)-b^{(1)}_t(y)|\leq    |x^{(1)}-y^{(1)}|^{\ff 2 3}\phi_n(|x^{(1)}-y^{(1)}|),\ \quad \text{if} \ x^{(2)}=y^{(2)},\\
&\|\nabla^{(2)}b^{(1)}_t(x)-\nabla^{(2)}b^{(1)}_t(y)\|\leq  \phi_n(|x^{(2)}-y^{(2)}|),\   \quad \text{if} \ x^{(1)}=y^{(1)}.\end{split}\end{equation*}
\item[{\bf(A3)}] {\bf (Regularity of $b^{(2)},\si$)}
Either
\beq\label{LB1'}\beg{cases}\beg{split}
 &|b^{(2)}_t(x)-b^{(2)}_t(y)| \leq  \big\{|x^{(1)}-y^{(1)}|^{\ff 2 3}\phi_n (|x^{(1)}-y^{(1)}|) +\phi_n^{\ff 72} (|x^{(2)}-y^{(2)}|)\big\},\\
 &\|\si_t(x)-\si_t(y)\|\le  |x-y| \ss{ \gg_n(|x-y|)}\,,\ \   |x|\lor |y|\le n; \end{split}\end{cases}\end{equation} or $\sup_{|x|\le n}\|\nn^{(2)}\si_t(x)\|_\infty\le C_n$ and for  $|x|\lor |y|\le n,$
\beq\label{LB2'}\beg{cases}\beg{split}
&|b^{(2)}_t(x)-b^{(2)}_t(y)| \leq  \big\{|x^{(1)}-y^{(1)}|^{\ff 2 3}\phi_n (|x^{(1)}-y^{(1)}|) +\phi_n  (|x^{(2)}-y^{(2)}|)\big\},\\
& \|\nn^{(2)} \si_t(x^{(1)},x^{(2)})-\nn^{(2)}\si_t(y^{(1)}, x^{(2)})\|\le   |x^{(1)}-y^{(1)}|\ss{\gg_n(|x^{(1)}-y^{(1)}|)}\,,\\
&\|\si_t(x^{(1)},x^{(2)})-\si_t(y^{(1)}, x^{(2)})\|\le  |x^{(1)}-y^{(1)}|\ss{\gg_n(|x^{(1)}-y^{(1)}|)}.
\end{split}\end{cases}\end{equation}
  \end{enumerate}

\beg{thm}\label{T1.1} \beg{enumerate}\item[$(1)$] Under assumption {\bf (A)},  for any $x\in \R^{d_1+d_2}$,
   SDE $\eqref{SDE}$ has a unique solution $ X_t(x) $ up to the explosion time $\zeta(x)$.
\item[$(2)$] If, in particular, $ b_t(x)$ and  $\si_t(x)$ do not depend on $x^{(1)}$, then the above  assertion follows provided for any $n\in\N$ there exists $\phi_n\in\D_0\cap\sS_0$ and $\gg_n\in \C$ such that
{\bf (A1)} and
\beq\label{LPT} \beg{cases} \|\si_t(x)-\si_t(y)\|\le |x^{(2)}-y^{(2)}|\ss{\gg_n(|x^{(2)}-y^{(2)}|)},\\ |b_t^{2}(x)-b_t^{(2)}(y)|+\|\nn^{(2)}b^{(1)}_t(x)
-\nn^{(2)}b^{(1)}_t(y)\|\le \phi_n(|x^{(2)}-y^{(2)}|) \end{cases} \end{equation} hold for all $t, |x|,|y|\le n$.
 \item[$(3)$] If there  exists $H\in C^2(\R^{d_1+d_2})$  such that
 \beq\label{1.7} H\ge 1,\ \ \lim_{|x|\to\infty} H(x)=\infty,\ \ |\nn^{(2)} H|^2\le C H^{2-\vv},\ \ \sL_t^{\Sigma,b} H\le \Phi(t) H,\ \ t\ge 0
 \end{equation}
 holds for some constant $\vv\in (0,1]$ and positive increasing function $\Phi$,
 then the solution to $\eqref{SDE}$ is non-explosive and for any $\eps'\in[0,\eps)$,
 \beq\label{1.8} \E\exp\bigg[\sup_{t\in [0,T]} H(X_t(x))^{\vv'} \bigg] \le
 \Psi(T) \exp\big[H(x)^\vv\big],\ \ T>0, x\in \R^{d_1+d_2}
 \end{equation}
 holds for some increasing function $\Psi: [0,\infty)\to (0,\infty).$
 \end{enumerate}  \end{thm}

\paragraph{Remark 1.2.}   (1) When $b^{(1)}$ is linear, an infinite-dimensional version of the well-posedness has been proved in \cite{Wa-Zh}   by following the line of  \cite{W15} for non-degenerate SPDEs, see \cite{AB1, AB2, AB3, AB4} for   discussions on the pathwise uniqueness of   SPDEs
with H\"older continuous drifts and non-degenerate additive noises.

(2) When $m=d$, the well-posedness was also  proved in  \cite{CDR} under a stronger assumption where  $\si$ is Lipschitz continuous,  $b(x)$ is H\"older continuous of order $\aa\in (\ff 3 2,1)$ in $x^{(1)}$  and  order $\bb\in (0,1)$ in $x^{(2)}$, and $\nn^{(2)} b^{(1)}$ is H\"older continuous.
In fact, we will show in Theorem \ref{T1.2} below that under this assumption and that $\si\in C^{1+\vv}$ for some $\vv>0$ the solutions to \eqref{SDE} form   $C^1$-stochastic diffeomorphism flows. Notice that the proofs given in  \cite{CDR} strongly depend on the explicit form of the fundamental solutions of linear degenerate Kolmogorov's operators,
%new
while our proof is based on  explicit probability   formulas of the semigroup associated to the linear stochastic Hamiltonian system (see Section 2.4 below).
%and new

\

To illustrate Theorem \ref{T1.1}, we present below three direct consequences, where the first generalizes to \eqref{SDE0}, the second includes a class of SDEs with unbounded time-delay which are interesting by themselves,  and the last presents a new well-posedness result for     non-degenerate SDEs.

\beg{cor}\label{C1.2} The following stochastic differential-integral equation on $\R^d$ admits a unique strong solution up to life time:
$$
\d X_t= \bigg(b_t(X_t)+\int_0^t\si_s(X_s)\d W_s\bigg)\d t,
$$
where $W_t$ is a $d$-dimensional Brownian motion,  $b: \R_+\times \R^d\to \R^d, \si: \R_+\times \R^d\to \R^d\otimes\R^d$ are measurable such that $b,\si$ and $\si^{-1}$ are locally bounded, and for any $n\ge 1$ there exist   $\phi_n\in \D_0\cap\sS_0$ and $\gg_n\in \C$ such that
for all $t,|x|, |y|\le n$,
  \beq\label{CD1} \beg{split}&|b_t(x)-b_t(y)|\le   |x-y|^{\ff 2 3}\phi_n(|x-y|),\\
 &|\si_t(x)-\si_t(y)|\le   |x-y|\ss{\gg_n(|x-y|)}.
 \end{split}\end{equation}
  \end{cor}
 \beg{proof} Let $\tt X_t^{(1)}=X_t, \tt X_t^{(2)}= \int_0^t\si(X_s)\d W_s.$ Then the equation reduces to \eqref{SDE} on $\R^{d+d}$ with
 $$\tt b_t^{(1)}(\tt x):=  b_t(\tt x^{(1)})+\tt x^{(2)}, \ \ \tt b_t^{(2)}:=0,\ \ \tt\si_t(\tt x):= \si_t(\tt x^{(1)}).$$ Obviously, the local boundedness of $b,\si$ and $\si^{-1}$ as well as
 \eqref{CD1} imply {\bf (A)}  with \eqref{LB1'} for $(\tt b,\tt \si)$.  Then the proof is finished by Theorem \ref{T1.1}(1).\end{proof}

 \beg{cor}\label{C1.2'} Let $b$ and $\si$ satisfy {\bf (A)} and let $b_t^{(1)}(x)= b_t^{(1)}(x^{(2)})$ not depend on $x^{(1)}$. Then for any $Y_0=y\in \R^{d_2}$, the following SDE with unbounded time-delay has a unique solution up to life time:
 $$\d Y_t= b_t^{(2)}\bigg(\int_0^t b_s^{(1)}(Y_s)\d s, Y_t\bigg)\d t +\si\bigg(\int_0^t b_s^{(1)}(Y_s)\d s, Y_t\bigg)\d W_t,\ \ Y_0=y.$$
 \end{cor}

 \beg{proof} Let $X_t^{(1)}=\int_0^t b_s^{(1)}(Y_s)\d s$ and $X_t^{(2)}= Y_t$. Then the SDE reduces to \eqref{SDE} with $X_0=(0,y)\in \R^{d_1+d_2}$. So, the desired assertion follows from Theorem \ref{T1.1}.\end{proof}

Finally, since existing well-posedness results for non-degenerated SDEs at least assumed   that $\si$ is weakly differentiable
(see \cite{Fe-Fl, Zh0} and references within), the following result is new even in the non-degenerate setting.

\beg{cor}\label{C1.3} The following SDE on $\R^d$ admits a unique strong solution up to life time:
$$
\d X_t=  b_t(X_t) +\si_t(X_t)\d W_t,
$$
where $W_t$ is a $d$-dimensional Brownian motion,
$b: \R_+\times \R^d\to \R^d, \si: \R_+\times \R^d\to \R^d\otimes\R^d$ are measurable such that $b,\si$ and $\si^{-1}$ are locally bounded,
and for any $n\ge 1$ there exist   $\phi_n\in \D_0\cap\sS_0$ and $\gg_n\in \C$ such that for all $t,|x|, |y|\le n$,
 \beq\label{CD33} |b_t(x)-b_t(y)|\le   \phi_n (|x-y|), \
 |\si_t(x)-\si_t(y)|\le C(t) |x-y| \ss{\gg_n(|x-y|)}.  \end{equation}  \end{cor}
 \beg{proof} Let $\tt X_t^{(1)}=\int_0^t X_s\d s, \tt X_t^{(2)}= X_t.$ Then the equation reduces to \eqref{SDE} on $\R^{d+d}$ with
 $$\tt b_t^{(1)}(\tt x):=  \tt x^{(2)}, \ \ \tt b_t^{(2)}(\tilde x)=b_t(\tt x^{(2)}),\ \ \tt\si_t(\tt x)= \si_t(\tt x^{(2)}).$$ Obviously, the local
  boundedness of $b,\si$ and $\si^{-1}$, together with
 \eqref{CD33}, implies that {\bf (A1)} and  \eqref{LPT} for $(\tt b,\tt \si)$.
 Then the proof is finished by Theorem \ref{T1.1}(2).\end{proof}

\subsection{Stability of   solutions with respect to coefficients}

About the continuous dependence of strong solutions with respect to the coefficients $(b,\sigma)$, we have
\bt \label{T1.5}
Let $(b^k,\sigma^k)_{k\in\mN_\infty}$ be a sequence of functions satisfying {\bf (A1)}, {\bf (A2)} and
\beq\label{LB9}\beg{cases}\beg{split}
 &|(b^{k})^{(2)}_t(x)-(b^{k})^{(2)}_t(y)| \leq  \big\{|x^{(1)}-y^{(1)}|^{\ff 2 3}\phi_n (|x^{(1)}-y^{(1)}|) +\phi_n^{\ff 72} (|x^{(2)}-y^{(2)}|)\big\},\\
 &\|\si^k_t(x)-\si^k_t(y)\|\le  C_n|x-y|,\ \   t\leq n, |x|\lor |y|\le n
 \end{split}\end{cases}\end{equation}
with the same localization constants $C_n$ and  $\phi_n\in \D_0\cap\sS_0$.  Assume that
$(b^k,\sigma^k)$ satisfies \eqref{1.7} with the same $H$ and $C$,
and for each $t,x$,
 $$
\lim_{k\to\infty} \|\sigma^k_t(x)-\sigma^\infty_t(x)\|+|b^k_t(x)-b^\infty_t(x)|=0.
 $$
 Let $X^k_t(x)$ be the unique solution of \eqref{SDE} corresponding to $(b^k,\sigma^k)$ for each $k\in\mN_\infty$.
 Then for each $\eps, T>0$ and $x\in\mR^d$,
 \begin{align}\label{STA}
\lim_{k\to\infty}\mP\bigg(\sup_{t\in[0,T]}|X^k_t(x)-X^\infty_t(x)|\geq\eps\bigg)=0.
 \end{align}
Moreover, if for some $p>d$ and for all $T,R>0$,
\begin{align}\label{STA3}
\sup_{k\in\mN_\infty}\sup_{|x|\leq R}\mE\bigg(\sup_{t\in[0,T]}|\nabla X^k_t(x)|^p\bigg)<\infty,
\end{align}
then for each $\eps, R, T>0$,
\begin{align}\label{STA4}
\lim_{k\to\infty}\mP\bigg(\sup_{t\in[0,T]}\sup_{|x|\leq R}|X^k_t(x)-X^\infty_t(x)|\geq\eps\bigg)=0.
\end{align}
\et

\paragraph{Remark 1.3.}
See Theorem \ref{T1.2} below  for sufficient conditions of \eqref{STA3}.
According to \cite[Theorem 2.3]{Xi-Zh}, condition \eqref{STA3} can be replaced with the following weaker one:
for some $p>d$ and for all $T,R>0$,
$$
\sup_{k\in\mN_\infty}\mE\bigg(\sup_{t\in[0,T]}|X^k_t(x)-X^k_t(y)|^p\bigg)\leq C|x-y|^p,\ \ |x|\vee|y|\leq R.
$$
\subsection{ $C^1$-stochastic diffeormorphism flow}

In order to show the $C^1$-diffeomorphism flow property of $X_t(x)$,
%we need  the following stronger  conditions.
we need stronger  conditions as shown in the following result.
\beg{thm}\label{T1.2} Assume {\bf (C1)} and that for some constant $\bb\in (0,\ff 1 3)$ and increasing function $C: [0,\infty)\to \R_+ $ the conditions
\beg{equation*}\beg{split}
 &|b^{(1)}_t(x)-b^{(1)}_t(y)|\leq  C(t) |x^{(1)}-y^{(1)}|^{\bb+\ff 2 3},\ \quad \text{if} \ x^{(2)}=y^{(2)},\\
&\|\nabla^{(2)}b^{(1)}_t(x)-\nabla^{(2)}b^{(1)}_t(y)\|\leq C(t) |x^{(2)}-y^{(2)}|^\bb,\   \quad \text{if} \ x^{(1)}=y^{(1)},\\
&|b^{(2)}_t(x)-b^{(2)}_t(y)| \leq C (t)\big( |x^{(1)}-y^{(1)}|^{\bb+\ff 2 3}  + |x^{(2)}-y^{(2)}|^\bb\big),\\
&\|\nn \si_t\|_\infty\le C(t),\ \ \|\nn\si_t(x)-\nn\si_t(y)\|\le C(t)|x-y|^{\bb}
\end{split}\end{equation*} hold for any $|x-y|\le 1, t\ge 0.$
Then the unique strong solution $\{X_t(\cdot)\}_{t\geq 0}$ to $\eqref{SDE}$  is a $C^1$-stochastic diffeomorphism flow, and
\beq\label{FL}
\sup_{x\in \R^{d_1+d_2}}  \E \bigg(\sup_{t\in [0,T]}  \|\nn X_t(x)\|^p\bigg) <\infty,\ \ T>0, p\ge 1.
\end{equation} \end{thm}

In the above result, $b$  has at most linear growth.
The following result shows that by making perturbations to $b$, it is possible to prove the $C^1$-stochastic diffeomorphism flow property for $b$ of high order polynomial growth.

\beg{thm}\label{T1.2'}
Keep the same assumptions of Theorem $\ref{T1.2}$.
Let $a:\mR_+\times\mR^{d_1+d_2}\to\mR^{d_2}$ be a measurable function such that $\nn a_t$ is locally H\"older continuous uniformly in $t\in[0,T]$ for any $T>0$.
Suppose also that for some $H\in C^2(\R^{d_1+d_2})$, $\vv\in (0,1]$, $\delta_1,\delta_2, C_1,C_2,C_3>0$ and positive increasing function $\Phi$,
and for all $t\geq 0$, $x\in\mR^{d_1+d_2}$,
 \beq\label{1.20}
 %H\ge 1,\ \lim_{|x|\to\infty} H(x)=\infty,\ |\nn^{(2)} H|^2\le C H^{2-\vv},\ |\nabla H|\leq CH,\ |\sL_t^{\Sigma,b+a} H|\le \Phi(t) H,
C_1(1+|x|^{\delta_1})\le H(x)\le C_2(1+|x|^{\delta_2}),\ |\nn^{(2)} H|^2\le C_3 H^{2-\vv},\ |\sL_t^{\Sigma,b+a} H|\le \Phi(t) H,
 \end{equation}
and for some $\eps'\in[0,\eps)$ and positive increasing function $\Phi'$, and for all $t\geq 0$ and $x, x'\in\mR^{d_1+d_2}$,
\begin{align}
|a_t(x)|\leq \Phi'(t) H(x)^{\eps'},\ |a_t(x)-a_t(x')|\leq \Phi'(t)(H(x)^{\eps'}+H(x')^{\eps'})|x-x'|,\label{NB4}
\end{align}
Then the  SDE
\begin{align}\label{SDE1}
\dif X_t=[a_t(X_t)+b_t(X_t)]\dif t+(0,\sigma_t(X_t)\dif W_t),\ \ X_0=x\in\R^{d_1+d_2}
\end{align}
has a unique strong solution   $X_t(x)$   such that $\{X_t(\cdot)\}_{t\geq 0}$ forms a $C^1$-stochastic diffeomorphism flow,
and for any $T>0$ and $p\ge 1$, there exists a constant $C>0$ such that
\beq\label{FL'}
\E\bigg(\sup_{t\in [0,T]} \|\nn X_t(x)\|^p\bigg) \le C\e^{H(x)^{\vv}},\ \ x\in \R^{d_1+d_2}.
\end{equation}
\end{thm}

Below is a simple example illustrating  Theorem \ref{T1.2'}, where   the drift   is neither local Lipschitz nor of linear growth.

\paragraph{Example 1.1.} Let $d_1=d_2=d$,  $\alpha\in(\frac{2}{3},1]$, $m\in\mN$ and $c_1,c_2>0.$ Take
$$
H(x)=1+\tfrac{ 1}{ 2} |x^{(2)}|^2+ c_1|x^{(1)}|^{\alpha+1}+c_2|x^{(1)}|^{m+1}.
$$
Let $\sigma$ be an invertible  $d\times d$-matrix. Consider the following SDE
$$
\dif (X^{(1)}_t, X^{(2)}_t)=(X_t^{(2)}, -\nn^{(1)}  H(X_t)) \dif t+(0,\sigma\dif W_t).
$$
It is easy to see that   Theorem \ref{T1.2'} applies to
$$
b(x)=\big(x^{(2)}, -c_1(\aa+1)x^{(1)}|x^{(1)}|^{\alpha-1}\big),\ \ a(x)=\big(0, -c_2(m+1)x^{(1)}|x^{(1)}|^{m-1}\big).
$$

\

In the spirit of \cite{Zv, Ve}, the key point of the study is to construct a time-dependent diffeomorphism on $\R^{d_1+d_2} $ which transforms   \eqref{SDE} into an equation with regular enough coefficients ensuring the desired assertions.   To this end, we take a freezing coefficient argument, which is different from the one used in \cite{CDR}, so that the construction is reduced to solve an parabolic equation associated to a linear stochastic Hamiltonian system.
To   figure out the minimal  conditions   on $b$ and $\si$ for the required estimates on solutions to this parabolic equation,
 we introduce some techniques in Section 2,   in particular, some  characterizations of the continuity using the heat semigroup. Moreover,
  in Section 2 we also present gradient estimates on the semigroup of the linear stochastic Hamiltonian system. With these preparations,
  in Section 3 we investigate   the   parabolic equation associated to the generator $\L_t^{\Sigma,b}$ (see \eqref{Eq1} below), which in turn   provides the desired   diffeomorphism on $\R^{d_1+d_2}$. Finally, in Section 4 we present complete proofs of the above theorems.

\section{Preparations}

This section contains some results which will be used to construct the regularization transform in the proof of the main results. We first present a Volterra-Gronwall  type inequality associated to a Dini function, then characterize the continuity of functions using the heat semigroup, and finally
introduce derivative formula and gradient estimates on linear stochastic Hamiltonian systems.

Throughout the paper, the letter $C$ with or without subscripts
will denote a positive constant whose value   may change from one appearance to another.
For two real functions $f$ and $g$, we write $f \preceq g $ if $f \leq C_0 g $ for some $C_0>0$; and $f \asymp g  $ if $C_1 g \leq f \leq C_2 g $
for some $C_1,C_2>0$.

 \subsection{Volterra-Gronwall  inequality associated to a Dini function}

\beg{lem} \label{Le23}
Let  $\phi: \R_+\to \R_+$ be a Dini function. For any $T>0$, there exists a constant $C=C(\phi, T)>0$ such that
if $\lambda\geq 0$ and bounded measurable functions $f, h:\R_+\to \R_+$ satisfy
$$
h(t)\leq \int^t_0\e^{-\lambda(t-s)}\ff{\phi(t-s)}{t-s} \big(h(s)+f(s)\big)\dif s, \ \ t\in (0,T],
$$
then
$$
h(t)\leq C\int^t_0\e^{-\lambda(t-s)}\ff{\phi(t-s)}{t-s}f(s)\dif s,\ \ t\in (0,T].
$$
\end{lem}

\begin{proof} Let $a_1(t)= \ff{\phi(t)}t$ and define
$$a_{n+1}(t)= \int_0^t a_n(t-s) a_1(s) \d s,\ \ t\in (0,T],\ n\in\mN.$$
Since $\int^T_0\ff{\phi(t)}t\dif t<\infty$, by   [27, Theorem 1] with $k(t,s):= \ff{\phi(t-s)}{t-s}1_{\{s<t\}}$ (see also \cite[Lemma 2.1]{Zh1}), we have
$$
a(t):=\sum_{n=1}^\infty a_n(t)\in L^1([0,T])
$$
and
\beq\label{Eb9} \beg{split}
a(t)=a_1(t)+\int^t_0a(t-s)a_1(s)\dif s.
\end{split}\end{equation}
Letting
$$
g(t)=\int^t_0\e^{-\lambda(t-s)}a_1(t-s)f(s)\dif s,
$$
then by \cite[Lemma 2.2]{Zh1}, we have
$$
h(t)\leq g(t)+\int^t_0 \e^{-\ll (t-s)} a(t-s)g(s)\dif s.
$$
Combining this with \eqref{Eb9} and using    Fubini's theorem, we obtain
\begin{equation*} \beg{split}
h(t)&\leq g(t)+\int^t_0 \e^{-\lambda (t-s)}a(t-s)\left(\int^s_0\e^{-\lambda(s-r)}a_1(s-r) f(r)\dif r\right)\dif s\\
&=g(t)+\int^t_0\left(\int^t_r \e^{-\lambda (t-s)} a(t-s)\e^{-\lambda(s-r)}a_1(s-r) \dif s\right)f(r)\dif r\\
&=g(t)+\int^t_0\e^{-\lambda (t-r)} f(r)\dif r\int_0^{t-r} a(t-r-s)a_1(s)\d s\\
&\le g(t)+\int^t_0\e^{-\lambda (t-r)}a(t-r) f(r)\dif r.
\end{split} \end{equation*} So, it remains to prove
\beq\label{FD}a(t)\le C a_1(t),\ \ t\in (0,T]\end{equation} for some constant $C>0.$
By the increasing property of $\phi$, we have
$$a_1(r t)=\ff{\phi(r t)}{r t} \le \ff{\phi(t)}{r t} =\ff{a_1(t)}r,\ r\in (0,1),\ t\in (0,T].$$
By the standard induction argument, this implies
\begin{align}a_n(r t) \le \ff{a_n(t)}r,\ \ r\in (0,1),\  t\in (0,T], \ n\in\mN.\label{YY1}\end{align}
Indeed, by the change of variables and induction hypothesis, we have
\begin{align*}
a_{n+1}(r t)&=\int_0^{r t} a_n(r t-s) a_1(s) \d s
=r\int_0^{t} a_n(r(t-s)) a_1(r s) \d s\\
&\leq\frac{1}{r}\int_0^{t} a_n(t-s) a_1(s) \d s=\frac{a_{n+1}(t)}{r}.
\end{align*}
Thus, for any $\eps\in(0,1)$ and $t\in(0,T]$, by \eqref{YY1} we have
\beg{equation*}\beg{split} &\int_0^t a(t-s)a_1(s)\d s = \sum_{n=1}^\infty \int_0^t a_n(t-s)a_1(s)\d s\\
&\le \sum_{n=1}^\infty \int_{\vv t}^t a_n(t-s)\ff{ a_1(t)}{s/t}\d s+ \sum_{n=1}^\infty \int_0^{\vv t} \ff{a_n(t)}{(t-s)/t}a_1(s)\d s\\
&\le \ff{a_1(t)}\vv \int_0^T a(s)\d s + \ff{a(t)}{1-\vv}\int_0^{\vv T} a_1(s)\d s.\end{split}\end{equation*}
Letting $\vv\in (0,1)$ be small enough such that $\ff{1}{1-\vv}\int_0^{\vv T} a_1(s)\d s\le \ff 1 2$, and combining this with \eqref{Eb9}, we obtain $$  a(t) \le 2 a_1(t)\bigg(1+ \ff 1 \vv \int_0^T a(t)\d t\bigg),\ \ t\in (0,T].$$ This implies \eqref{FD} since $a\in L^1([0,T]).$
\end{proof}

\subsection{Slowly varying functions}

We first recall some important properties of slowly varying functions (cf. \cite[Theorem 1.5.6 (ii) and Theorem 1.5.11]{Bi-Go-Te}).

\bp\label{Pro1}
For any $\phi\in\sS_0$, the following assertions hold:
\begin{enumerate}[{\rm (i)}]
\item For any $\delta>0$, there is a constant $C=C(\delta)\geq 1$ such that for all $t,s>0$,
\begin{align*}
\frac{\phi(t)}{\phi(s)}\leq C\max\left\{\Big(\frac{t}{s}\Big)^{\delta}, \Big(\frac{t}{s}\Big)^{-\delta}\right\}.
\end{align*}
\item For any $\beta>-1$, as $t\to 0$, we have
\begin{align*}
\int^t_0 s^{\beta}\phi(s)\dif s\sim\frac{t^{\beta+1}\phi(t)}{\beta+1},\ \ \int^1_ts^{-\beta-2}\phi(s)\dif s\sim \frac{t^{-\beta-1}\phi(t)}{\beta+1}.
\end{align*}
%\item If $\int^t_0s^{-1}\phi(s)\dif s<\infty$, then the function $t\mapsto \int^t_0s^{-1}\phi(s)\dif s$ is slowly varying.
\end{enumerate}
\ep

The following lemma is simple.
\bl\label{Le24}
For any bounded measurable function $\psi:(0,1]\to\R_+$, we have
\begin{align}\label{DB8}
[f]_\psi:=\sup_{|x-y|\leq 1}\frac{|f(x)-f(y)|}{\psi(|x-y|)}=\sup_{x\not= y}\frac{|f(x)-f(y)|}{  \psi_{[0]}(|x-y|)},
\end{align}
where $\psi_{[0]}(t):=\psi(t)1_{t\leq 1}+\psi_*(1)t1_{t>1}$ and $\psi_*(1):=\sup_{s\in(0,1]}\psi(s)$.
\el
\begin{proof}
Clearly, it suffices to prove that
$$
|f(x)-f(y)|\leq [f]_\psi\psi_*(1)|x-y|, \ |x-y|\geq 1.
$$
Suppose that $n<|x-y|\leq n+1$ for some $n\in\mN$. Let $x=x_0, x_1,\cdots, x_{n}, x_{n+1}=y$ be $n+2$-points in $\mR^d$ so that
$$
|x_i-x_{i-1}|=1,\  i=1,\cdots,n,\  |x-y|=n+|x_{n+1}-x_n|.
$$
Then we have
$$
|f(x)-f(y)|\leq \sum_{i=1}^{n+1}|f(x_i)-f(x_{i-1})|\leq [f]_\psi\psi_*(1)(n+|x_{n+1}-x_n|)=[f]_\psi\psi_*(1)|x-y|.
$$
The proof is finished.
\end{proof}

Due to the above lemma and also for later use, we introduce
\begin{align}\label{DD1}
\sR_\alpha:=\Big\{\phi_{[\aa]}(t):=t^\alpha\phi(t)1_{t\leq 1}+c_\alpha t1_{t>1}:  \phi\in\sS_0\ \mathrm{with}\ c_\alpha= \mathrm{sup}_{s\in(0,1]}(s^\alpha\phi(s))<\infty\Big\}
\end{align}
for $\alpha\in[0,1]$, and let
$$
\sR =\cup_{\alpha\in[0,1]}\sR_\alpha.
$$
The function $\phi_{[\aa]}$ with $\aa\in [0,1]$ and $\phi\in \D_0$ not only   characterizes the H\"older-Dini modulus, but also reduces the study to functions with linear growth. Notice that by (i) of Proposition \ref{Pro1}, $c_\alpha$ in \eqref{DD1} is automatically finite for $\alpha\in(0,1]$.

Below we list the main properties of $\psi\in\sR_\alpha$ for later use,
which are easy consequences of Proposition \ref{Pro1}.
\bp\label{Pro2}
For $\alpha\in[0,1]$, let $\psi\in\sR_\alpha$.
\begin{enumerate}[{\rm (i)}]
\item For any $\delta>0$, there is a constant $C=C(\delta)\geq 1$ such that for all $t,s>0$,
\begin{align}
\frac{\psi(t)}{\psi(s)}\leq C\max\left\{\Big(\frac{t}{s}\Big)^{\alpha+\delta}, \Big(\frac{t}{s}\Big)^{\alpha-\delta}\right\}.\label{Pr}
\end{align}
In particular, if $\alpha\in[0,1)$, then for all $t\geq s>0$,
\begin{align}\label{Pr0}
\frac{s}{\psi(s)}\leq C\frac{t}{\psi(t)}.
\end{align}
\item If $\alpha\in(0,1)$, then there is a constant $C>0$ such that for all $t\in(0,1]$,
\begin{align}
\int^t_0 s^{-1}\psi(s)\dif s\leq C\psi(t),\ \ \int^1_ts^{-2}\psi(s)\dif s\leq Ct^{-1}\psi(t).\label{Pr1}
\end{align}
%\item If $\alpha=0$ and $\int^1_0s^{-1}\psi(s)\dif s<\infty$, then $t\mapsto \tilde\psi(t):=1_{t\leq 1}\int^{t}_0s^{-1}\psi(s)\dif s+t1_{t\geq 1}\int^{1}_0s^{-1}\psi(s)\dif s$
%belongs to $\sR_0$.
\item There is a constant $C>0$ such that for all $s,t>0$,
\beg{equation}\label{JJ}\beg{split}
\psi(s+t)  \le C\big(\psi(s)+\psi(t)\big).
\end{split}\end{equation}
\end{enumerate}
\ep

\subsection{Characterization of continuity by using heat semigroup}

Let $\B_p(\R^d)$ be the set of all measurable functions on $\R^d$ with polynomial growth.
We will investigate the continuity of  $f\in \B_p(\R^d)$ on $\R^d$ by using the standard  heat semigroup
 \begin{align}
\bP_\theta f(x)=\int_{\mR^d}f(y)p_\theta(x-y)\dif y,\ \theta>0,\label{PT}
\end{align}
where
$$
p_\theta(x):=\ff 1 {(2\pi \theta)^{d/2}} \e^{-\ff{|x|^2}{2\theta}}.
$$
Notice that by elementary calculus,
\begin{align}
|\nabla^k\p^j_\theta p_\theta(x)|\preceq \frac{|x|^k(\theta+|x|^2)^j}{\theta^{k+2j}}p_\theta(x),\ \theta>0,\ x\in\mR^d,\ k,j=0,1.\label{LM1}
\end{align}
For any measurable function $\psi:[0,1]\to\R_+$ and $f:\mR^d\to\mR$, define
$$
[f]_{\psi}:=\sup_{|x-y|\leq 1}\frac{|f(x)-f(y)|}{\psi(|x-y|)},\  \|f\|_\infty:=\sup_{x\in\mR^d}|f(x)|,\ \|f\|_\psi:= [f]_\psi+\|f\|_\infty.
$$
It should be noticed by \eqref{DB8} and \eqref{DD1}  that for any $\psi\in\sR$,
\begin{align}\label{DB88}
|f(x)-f(y)|\leq \psi(|x-y|)[f]_\psi,\ \ x,y\in\mR^d,
\end{align}
and if $\psi_1(s)\leq C\psi_2(s), s\in(0,1]$ for some $C>0$, then
$$
[f]_{\psi_2}\leq C[f]_{\psi_1}.
$$

We first present the following simple lemma.

\beg{lem}\label{L0} For any $\psi\in\sR$ and $\beta\geq 0$, there exists a constant $C>0$ such that for all $\theta>0$,
\beg{align}  &\qquad\int_{\R^d}|z|^\beta \psi(|z|) p_\theta (z)\d z\le C \theta^{\frac{\beta}{2}}\psi(\theta^{\ff 1 2}),\label{Ed1}\\
&\|\nn^k\pp^{j}_\theta  \bP_\theta f\|_\infty \le C[f]_\psi \theta^{-\ff k 2-j}\psi(\theta^{\ff 1 2}), \ k, j=0,1.\label{Ed1'}
\end{align}
\end{lem}

\beg{proof}
 Let $\psi\in\sR_\alpha$ for some $\alpha\in[0,1]$.
 By the change of variables and \eqref{Pr}, for any $\dd\in(0,1)$, we have
\begin{align*}
&\int_{\R^d} |z|^\beta\psi(|z|)  p_\theta (z)\d z
=\theta^{\frac{\beta}{2}}\int_{\R^d} |z|^j\psi (\theta^{\ff 1 2}|z|)  p_1 (z)\d z\\
&\qquad\preceq \theta^{\frac{\beta}{2}}\psi (\theta^{\ff 1 2})\int_{\R^d}|z|^\beta\Big(|z|^{\alpha+\delta}\vee|z|^{\alpha-\delta}\Big) p_1 (z)\d z,
\end{align*}
which gives \eqref{Ed1}.

Next, for any $x\in \R^d$, let $f_x= f-f(x)$. By \eqref{LM1}, \eqref{DB88} and \eqref{Ed1} we obtain
\beg{align*} |\nn^k \p^j_\theta\bP_\theta f|(x)&= |\nn^k \p^j_\theta\bP_\theta f_x|(x)\le \int_{\R^d} |f_x(x+z)||\nn^k\p^j_\theta p_\theta (z)|\d z\\
 &\ppp [f]_\psi \int_{\R^d}\frac{|z|^k(\theta+|z|^2)^j\psi(|z|)}{\theta^{k+2j}}  p_\theta (z)\d z\ppp [f]_\psi\theta^{-\ff k 2-j}\psi(\theta^{\ff 1 2}).\end{align*}
 This proves \eqref{Ed1'}.
\end{proof}

We have the following commutator estimate result. A similar version for the Cauchy semigroup can be found in \cite{Ch-So-Zh}.
As an   advantage of the present result,   it applies to $f\in \B_p(\R^d)$, the class of measurable functions with polynomial growth.

\bl\label{Le31} Let $\psi\in\sR_\alpha$ for some $\alpha\in[0,1]$ and $\phi: \R_+ \to\R_+ $ be increasing so that $\psi\phi$ satisfies for some $C>0$,
\begin{align}\label{JJ0}
(\psi\phi)(t+s)\leq C\big((\psi\phi)(t)+(\psi\phi)(s)\big),\ \ t,s>0.
\end{align}
Suppose also that $\psi(t)$ is increasing on $[0,1]$ if $\alpha=0$, and $t^{-1}\psi(t)$ is decreasing on $[0,1]$ if $\alpha=1$.
Then there exists a  constant $C>0$ such that for any $f \in\sB_p(\mR^d)$ and $g\in \B_b(\R^d)$,
\beq\label{CC1}
[\pp_\theta\bP_\theta(fg)-f\pp_\theta\bP_\theta g]_{\psi}
\le C[f]_{\psi\phi}\|g\|_\infty  \theta^{-1}\phi(\theta^{\ff 1 2}),\ \ \theta\in(0,1].
\end{equation}
\el

\beg{proof}
By definition \eqref{PT}, we have
\beq\label{P1}
F_\theta(x):= \pp_\theta\bP_\theta (fg)(x)- f(x) \pp_\theta\bP_\theta g(x)= \int_{\R^d} (f(z)-f(x)) g(z) \pp_\theta p_\theta (x-z)\d z,
\end{equation}
which, by \eqref{DB88}, \eqref{LM1} and \eqref{Ed1}, implies that for all $\theta>0$,
\beg{equation}\label{P0}\beg{split}
\|F_\theta\|_\infty&\ppp [f]_{\psi\phi}
\|g\|_\infty\int_{\R^d}(\psi\phi)(|x-z|)|\pp_\theta p_\theta(z-x) |\,\d z\ppp [f]_{\psi\phi} \|g\|_\infty\theta^{-1}(\psi\phi)(\theta^{\ff 12}).
\end{split}\end{equation}
Thus, when $1\geq |x-y|^2\ge\theta$,  by \eqref{Pr} for $\alpha\in(0,1]$ and by the increasing property of $\psi$ for $\alpha=0$, we have
\beq\label{P2'}\beg{split}  |F_\theta(x)-F_\theta(y)|&\le 2\|F_\theta\|_\infty \ppp
[f]_{\psi\phi} \|g\|_\infty   \psi(|x-y|)\theta^{-1}\phi(\theta^{\frac{1}{2}}).
\end{split}\end{equation}
On the other hand, by \eqref{P1} we have
\beq\label{P3}\beg{split} F_\theta(x)-F_\theta(y)
&= \int_{\R^d} (f(z)-f(x)) g(z) (\pp_\theta p_\theta(x-z) -\pp_\theta p_\theta(y-z))\d z\\
&+\int_{\R^d} (f(y)-f(x)) g(z) \pp_\theta p_\theta (y-z)\d z=: I_1+I_2.
\end{split}\end{equation}
When $|x-y|^2\le\theta\leq 1$, by \eqref{DB88}, \eqref{JJ0}, \eqref{LM1} and \eqref{Ed1}, we have
\beq\label{P4} \beg{split}|I_1|& \ppp [f]_{\psi\phi} \|g\|_\infty |x-y| \int_{\R^d\times [0,1]} (\psi\phi)(|x-z|) |\nn \pp_\theta p_\theta(x-z+ r(y-x))|\d z\d r\\
&\ppp [f]_{\psi\phi} \|g\|_\infty|x-y| \int_{\R^d\times [0,1]} \Big[(\psi\phi)(|x-z+ r(y-x)|)+(\psi\phi)(|x-y|)\Big]\\
&\qquad \qquad\qquad \times  |\nn \pp_\theta p_\theta(x-z+ r(y-x))|\,\d z\d r\\
&\ppp [f]_{\psi\phi} \|g\|_\infty |x-y|(\psi\phi)(\theta^{\ff 1 2})\theta^{-\ff {3} 2}
\ppp [f]_{\psi\phi} \|g\|_\infty\psi(|x-y|)\theta^{-1}\phi(\theta^{\frac{1}{2}}),
  \end{split}\end{equation}
  where the last step is due to \eqref{Pr0} for $\alpha\in[0,1)$ and the decreasing property of $t^{-1}\psi(t)$ for $\alpha=1$.
Moreover, since $\phi$ is increasing,   when $|x-y|^2\le \theta$, it follows from \eqref{DB88},  \eqref{Ed1} that
   $$|I_2| \ppp [f]_{\psi\phi} \|g\|_\infty (\psi\phi)(|x-y|)  \theta^{-1}\leq [f]_{\psi\phi} \|g\|_\infty\psi(|x-y|)
\theta^{-1}\phi(\theta^{\ff 12}).
   $$
Combining this with \eqref{P2'}, \eqref{P3} and \eqref{P4}, we obtain \eqref{CC1}.
\end{proof}

We are now able to characterize a H\"older-Dini continuous function by using the heat semigroup (see \cite{St} for the characterization of H\"older space by using Poisson integrals).
\beg{lem}\label{LN} For any $\phi\in\sR$ with $\int^1_0\ff{\phi(s)}s\d s<\infty$, letting
\begin{align}\label{KL1}
\bar\phi(t)=t+t\int^1_t \frac{\phi(s)}{s^2}\dif s+\int_0^t\ff{\phi(s)}s\d s,\ t\in(0,1),
\end{align}
then we have
\beq\label{*WF}
\|f\|_{\bar\phi}\ppp \|f\|_\infty +  \sup_{\theta\in(0,1]} \left(\frac{\|\theta \pp_\theta\bP_\theta f\|_\infty}{\phi(\theta^{\frac{1}{2}})}\right),\ \ f\in\sB_b(\R^d).
\end{equation}
In particular, if $\phi\in \sR_\aa$ for some $\aa\in (0,1)$, then
\begin{align}\label{KL2}
\|f\|_{\phi}\asymp \|f\|_\infty +  \sup_{\theta\in(0,1]} \left(\frac{\|\theta \pp_\theta\bP_\theta f\|_\infty}{\phi(\theta^{\frac{1}{2}})}\right),\ \ f\in\sB_b(\R^d).
\end{align}
\end{lem}
\beg{proof}
Notice that
$$
f(x)=\bP_\theta f(x)-\int_0^\theta\pp_s\bP_sf(x)\d s.
$$
Since $\|\nabla\bP_s f\|_\infty\preceq\|f\|_\infty/\sqrt{s}$ for $s>0$ and $\p_s\nabla\bP_s f(x)=\nabla\bP_{s/2}(\p_r \bP_r)_{r=s/2}f(x)$,
we have
$$
\nabla\bP_\theta f(x)=\int^\infty_\theta \p_s\nabla\bP_s f(x)\dif s=\int^\infty_\theta \nabla\bP_{s/2}(\p_r \bP_r)_{r=s/2}f(x)\dif s,
$$
which, by \eqref{Ed1'},  implies that for $\theta\in(0,1]$,
$$
\|\nabla\bP_\theta f\|_\infty\preceq \ell(f)\left(\int^\infty_1 s^{-\ff 3 2}\dif s+\int^1_\theta s^{-\ff 3 2}\phi(s^{\ff 1 2})\dif s\right)
\preceq \ell(f)\left(1+\int^1_{\ss\theta} s^{-2}\phi(s)\dif s\right),
$$
where $\ell(f)$ is the quantity of the right hand side of \eqref{*WF}. Hence,
\begin{align*}
|f(x)-f(y)|&\leq \|\nabla\bP_\theta f\|_\infty|x-y|+2\int_0^\theta\|\pp_s\bP_sf\|_\infty\d s\\
&\preceq \ell(f)\left(|x-y|+|x-y|\int^1_{\ss\theta} s^{-2}\phi(s)\dif s+\int_0^\theta s^{-1}\phi(s^{\ff 1 2})\d s\right),
\end{align*}
which in turn implies that by letting $\theta=|x-y|^2\leq 1$,
$$
|f(x)-f(y)|\preceq \ell(f)\bar\phi(|x-y|),
$$
where $\bar\phi$ is defined by \eqref{KL1}. If $\alpha\in(0,1)$, by \eqref{Pr0} and \eqref{Pr1},
we have $\bar\phi(t)\preceq\phi(t)$. Thus, \eqref{KL2} follows by \eqref{*WF} and \eqref{P0}
with $g=1$ and $\psi=1$.
 \end{proof}

Next, we consider the product space $\R^{d_1+d_2}$. For any $\psi_1,\psi_2: \R_+\to \R_+$ and $f\in C(\R^{d_1+d_2})$, set
\beg{equation*}\beg{split}
[f]_{\psi_1,\infty}:= \sup_{x^{(2)}\in\R^{d_2}} [f(\cdot,x^{(2)})]_{\psi_1},\ \ &\  [f]_{\infty,\psi_2}:= \sup_{x^{(1)}\in\R^{d_1}} [f(x^{(1)},\cdot)]_{\psi_2},\\
[f]_{\psi_1,\psi_2}:= [f]_{\psi_1,\infty}+ [f]_{\infty,\psi_2},\ \ & \ \|f\|_{\psi_1,\psi_2}:= [f]_{\psi_1,\psi_2}+ \|f\|_\infty,
\end{split}\end{equation*}
and for simplicity,
$$
[f]_{\psi}:= [f]_{\psi,\psi},\ \ \|f\|_{\psi}:= \|f\|_{\psi,\psi}.
$$
Let $\bP_\theta^{(i)}$ be the heat semigroup on $\R^{d_i}$, we set for $x=(x^{(1)},x^{(2)})\in \R^{d_1+d_2}$,
\beq\label{PP11}
\bP_\theta^{(1)} f(x)= \big\{\bP_\theta^{(1)} f(\cdot, x^{(2)})\big\}(x^{(1)}),\ \
\bP_\theta^{(2)} f(x)= \big\{\bP_\theta^{(2)} f(x^{(1)},\cdot) \big\}(x^{(2)}).
\end{equation}
Obviously, Lemmas \ref{Le31} and \ref{LN} apply to both $(\|\cdot\|_{\phi,\infty}, \bP_\theta^{(1)})$ and $(\|\cdot\|_{\infty, \phi}, \bP_\theta^{(2)}).$
For instance, letting $\bP_\theta=\bP_\theta^{(1)}\bP_\theta^{(2)}$ be the Gaussian heat semigroup on $\mR^{d_1+d_2}$,
by the contractivity of $\bP_\theta^{(i)}$ under the uniform norm, Lemma \ref{Le31} implies the following result.

\bl\label{Le26}
Let $\psi_1,\psi_2\in\sR$ and $\phi: \R_+ \to\R_+ $ be increasing such that $\psi_i, i=1,2$ and $\phi$ satisfy the same assumptions as in Lemma \ref{Le31}.
Then there exists a  constant $C>0$ such that
\beq\label{CC0}
[\pp_\theta\bP_\theta(fg)-f\pp_\theta\bP_\theta g]_{\psi_1,\psi_2}
\le C[f]_{\psi_1\phi,\psi_2\phi}\|g\|_\infty  \theta^{-1}\phi(\theta^{\ff 1 2}),\ \ \theta>0.
\end{equation}
\el

\

Finally, the following result characterizes $\|\cdot\|_{\infty,\phi}$ by using $\bP^{(1)}$, and the same holds for $(\|\cdot\|_{\phi,\infty}, \bP_\theta^{(2)}).$

\beg{lem}\label{LN9} For any $\psi_1,\psi_2\in \sR$ and $\phi:\R_+\to\R_+$, there exists a constant $C>0$ such that
\beg{align*}
&[\pp_\theta \bP_\theta^{(1)} (fg)- f\pp_\theta \bP_\theta^{(1)} g]_{\infty, \phi}\le C [f]_{\psi_1,\psi_2}\|g\|_{\infty,\psi_2}
 \sup_{s\in (0,1]}\bigg(\ff{\psi_1(\theta^{\ff 1 2})\psi_2(s)+\psi_1(\theta^{\ff 1 2})\wedge\psi_2(s)}{\phi(s)}\bigg)\theta^{-1}
\end{align*}
holds for all $\theta \in (0,1]$ and   measurable functions $f,g$ on $\R^{d_1+d_2}$.
\end{lem}
\beg{proof} By definition, we have
$$
F_\theta(x)= \pp_\theta \bP_\theta^{(1)} (fg)(x)- f(x)\pp_\theta \bP_\theta^{(1)} g(x)
= \int_{\mR^{d_1}}G(z^{(1)}, x^{(1)}, x^{(2)})\p_\theta p_\theta(x^{(1)}-z^{(1)})\dif z^{(1)},
$$
where
$$
G(z^{(1)}, x^{(1)}, x^{(2)}):=\Big(f(z^{(1)}, x^{(2)})-f(x^{(1)}, x^{(2)})\Big)g(z^{(1)}, x^{(2)}).
$$
Clearly, by \eqref{DB88} we have
\begin{align*}
&|G(z^{(1)}, x^{(1)}, x^{(2)})-G(z^{(1)}, x^{(1)}, y^{(2)})|\leq\psi_1(|x^{(1)}-z^{(1)}|)[f]_{\psi_1,\infty}\psi_2(|x^{(2)}-y^{(2)}|)[g]_{\infty,\psi_2}\\
&\qquad\qquad+2\Big(\big(\psi_1(|x^{(1)}-z^{(1)}|)[f]_{\psi_1,\infty}\big)\wedge\big(\psi_2(|x^{(2)}-y^{(2)}|)[f]_{\infty,\psi_2}\big)\Big)\|g\|_\infty.
\end{align*}
Hence, for $x,y\in\mR^{d_1+d_2}$ with $x^{(1)}=y^{(1)}$, by \eqref{Ed1}, we obtain
\begin{align*}
|F_\theta(x)-F_\theta(y)|&\preceq [f]_{\psi_1,\infty}[g]_{\infty,\psi_2}\psi_1(\theta^{\ff 1 2})\psi_2(|x^{(2)}-y^{(2)}|)\theta^{-1}\\
&+[f]_{\psi_1,\psi_2}\|g\|_{\infty}\Big(\psi_1(\theta^{\ff 1 2})\wedge\psi_2(|x^{(2)}-y^{(2)}|)\Big)\theta^{-1},
\end{align*}
which in turn gives the desired estimate by dividing both sides by $\phi(|x^{(2)}-y^{(2)}|)$ and then taking supremum for $|x^{(2)}-y^{(2)}|\leq 1$.
\end{proof}
\subsection{Gradient  estimates for linear stochastic Hamiltonian system}

Let $B: \R_+\to \R^{d_1}\otimes\R^{d_2},\sigma: \R_+\to\R^{d_2}\otimes\R^{d_2}$ be measurable such that $B_rB_r^*$ and $\si_r$ are invertible with
\beq\label{IV}
\kappa:=\sup_{r\in\R_+} \big(|B_r|+|\si_r|+|(B_rB_r^*)^{-1}|+|\si_r^{-1}|\big)<\infty.
\end{equation}
For $x=(x^{(1)},x^{(2)})\in\mR^{d_1+d_2}$ and $0\leq s\leq t$, define
\beq\label{SL}\beg{split}
X_{s,t}(x)=\left(x^{(1)}+\GG_{s,t}x^{(2)}+\int^t_s B_r\d r\int^r_s\sigma_{r'}\dif W_{r'}, x^{(2)}+\int^t_s\sigma_{r}\dif W_{r}\right),\end{split}
\end{equation}
where $(W_r)_{r\geq 0}$ is a $d_2$-dimensional standard Brownian motion, and
\begin{align}\label{GG}
\GG_{s,t} = \int_s^t B_r\dif r.
\end{align}
Clearly, $X_{s,t}(x)=(X^{(1)}_{s,t}, X^{(2)}_{s,t})$ solves the following degenerate linear equation for $t\ge s$:
\begin{align}
\left\{
\begin{aligned}
\dif X^{(1)}_{s,t}&=B_tX^{(2)}_{s,t}\dif t,& X^{(1)}_{s,s}=x^{(1)},\\
\dif X^{(2)}_{s,t}&=\sigma_t\dif W_t, &X^{(2)}_{s,s}=x^{(2)}.
\end{aligned}
\right.
\end{align}
Let $P_{s,t}$ be the Markov operator associated with $X_{s,t}(x)$, i.e.,
$$
P_{s,t} f(x)=\mE f(X_{s,t}(x)),\ f\in \B_p(\mR^{d_1+d_2}).
$$

We first investigate the derivative estimates of $P_{s,t}f.$ To this end, we collect some frequently used notations here.
\begin{itemize}
\item For a smooth function $f$ on $\mR^{d_1+d_2}$, $\nabla^{(1)}f$ and $\nabla^{(2)}f$ denotes the gradient of $f$
with respect to the variables $x^{(1)}$ and $x^{(2)}$ respectively. In particular, by \eqref{SL} we have
\beq\label{AA5}
 \nn^{(1)}P_{s,t} f=P_{s,t}\nn^{(1)}f,\ \   P_{s,t}\nn^{(2)} f= \nn^{(2)}P_{s,t} f- \GG_{s,t} \nn^{(1)}P_{s,t}f.
\end{equation}
\item For $h=(h^{(1)}, h^{(2)})\in\mR^{d_1+d_2}$, we also write
$$
\nabla f:=(\nabla^{(1)}f,\nabla^{(2)}f),\ \nabla_h f:=\<\nabla f,h\>=\nabla^{(1)}_{h^{(1)}} f+\nabla^{(2)}_{h^{(2)}} f.
$$
\item Let $\sU$ be the set of all increasing functions $\phi:\R_+\to\R_+$ with the property
\begin{align}
\phi( rt)\leq Cr^\delta\phi(t),\ \ r\geq 1, \ \ t>0\label{UU0}
\end{align}
for some $C,\delta>0$. Notice that by \eqref{Pr},
$$
\sD_0\cap\sR\subset\sU.
$$
\end{itemize}

To estimate the derivatives of $P_{s,t}f$, we first present a Bismut type derivative formula which can be found in \cite{Zh1},
\cite{Gu-Wa} and \cite{Wa-Zh}. For readers' convenience we state the formula in details and present a simple proof.

Fix $0\leq s\leq t$ and define
$$
Q_{s,t} = \int^t_s(t-r)(r-s)B_rB_r^*\dif r\in\R^{d_1}\otimes\R^{d_1}.
$$
By \eqref{IV}, it holds that for some $C>0$,
$$
|Q^{-1}_{s,t}|\leq C(t-s)^{-3},\ \ t>s.
$$
For $h=(h^{(1)},h^{(2)})\in\mR^{d_1+d_2}$, define for $r\in[s,t]$,
\begin{align}
\Phi^{h}_{s,t}(r)=\frac{h^{(2)}}{t-s}+(t+s-2r)B_s^*Q^{-1}_{s,t}\left[h^{(1)}+ \int^t_s\frac{t-r'}{t-s}B_{r'}h^{(2)} \dif r' \right].\label{Phi}
\end{align}
Obviously, by \eqref{IV}, there exists a constant $C>0$ such that
\beq\label{EST} \big| \Phi^{h}_{s,t}(r)\big|\le C \left(\ff{|h^{(2)}|}{t-s} + \ff{|h^{(1)}|}{(t-s)^2}\right),\ \ 0\le s<t<\infty, r\in [s,t].\end{equation}

\bt\label{T2.1}
For $n\in\mN$, $s=s_0<s_1\cdots<s_n=t$ and $h_1,\cdots, h_n\in \R^{d_1+d_2}$, let
\begin{align}\label{xi}
\breve h_{i}=\left(h_i^{(1)}+\GG_{s,s_{i-1}} h_i^{(2)},\ h^{(2)}_i\right),\ \
\xi_{s_{i-1},s_{i}}^{\breve h_i}=\int_{s_{i-1}}^{s_{i}}\<\sigma_r^{-1} \Phi^{\breve h_i}_{s,t}(r),\dif W_r\>,
\end{align}
where $\Gamma_{s, s_{i-1}}$ is defined by \eqref{GG} and $i=1,\cdots, n$.
Then for any $f\in \B_p(\mR^{d_1+d_2})$, we have
\begin{align}\label{BS}
\nabla_{h_1}\cdots\nn_{h_n}  P_{s,t} f(x) = \mE\left[ f\big(X_{s,t}(x)\big)\prod_{i=1}^{n}\xi_{s_{i-1},s_{i}}^{\breve h_{i}} \right],\ x\in\mR^{d_1+d_2}.
\end{align}
\et

\begin{proof}
(i) First of all, we consider the case of $n=1$. For $\eps\in(0,1)$, define
$$
W_r^\eps= W_r-\eps \int_{s}^r \sigma_{r'}^{-1} \Phi^{h}_{s,t}(r') \dif r',\ \ r\in [s,t].
$$
By Camaron-Martin's theorem, $(W^\eps_r)_{r\in[s,t]}$ is still a Brownian motion under the probability measure $\dif \mP_\eps:= R_\eps\dif\mP$, where
\begin{align}\label{2.7}
R_\eps:= \exp\left[\eps \int_{s}^{t} \<\sigma_r^{-1} \Phi^{h}_{s,t}(r), \dif W_r\>
-\frac{\eps^2} 2  \int_{s}^{t} \big|\sigma^{-1}_r\Phi^{h}_{s,t}(r)\big|^2\dif r\right].
\end{align}
Thus, if we write
$$
X^\eps_{s,t}(x):=\left(x^{(1)}+\eps h^{(1)}+\int^t_s B_r\left[x^{(2)}+\eps h^{(2)}+\int^r_s\sigma_{r'}\dif W^\eps_{r'}\right]\dif r,
x^{(2)}+\eps h^{(2)}+\int^t_s\sigma_{r}\dif W^\eps_{r}\right),
$$
then the law of $X_{s,t}(x+\eps h)$ under $\mP$ is the same as the law of $X^\eps_{s,t}(x)$ under $\mP_\eps$, that is,
$$
P_{s,t}f(x+\eps h)=\mE f(X_{s,t}(x+\eps h))=\mE(R_\eps f(X^\eps_{s,t}(x))).
$$
On the other hand, by definition \eqref{Phi}, it is easy to see that
\begin{align*}
X^\eps_{s,t}(x)=X_{s,t}(x)+\eps\left(h^{(1)}+\int^t_s B_r\left[h^{(2)}-\int^r_s\Phi^{h}_{s,t}(r')\dif r'\right]\dif r, h^{(2)}-\int^t_s\Phi^{h}_{s,t}(r)\dif r\right)=X_{s,t}(x).
\end{align*}
Hence,
$$
\nabla_h P_{s,t}f(x)=\lim_{\eps\downarrow 0} \frac{1}{\eps}\mE\left[f(X_{s,t}(x+\eps h))-f(X_{s,t}(x))\right]
=\lim_{\eps\downarrow 0} \mE\left[\frac{R_\eps-1}\eps f(X_{s,t}(x))\right],
$$
which together with \eqref{2.7} yields \eqref{BS} for $n=1$.
\\
\\
(ii) Assuming that \eqref{BS} holds for $n=k\in\mN$, we intend to prove \eqref{BS} for $n=k+1.$  Noticing that $P_{s,t} f=P_{s,s_k} P_{s_{k},t} f$ and by definition \eqref{SL},
\begin{align*}
\nabla_{h_{k+1}}X_{s,s_k}=\left(h_{k+1}^{(1)}+\GG_{s,s_k} h_{k+1}^{(2)},\ h^{(2)}_{k+1}\right)=\breve h_{k+1},
\end{align*}
by  induction hypothesis, we have
\begin{align*}
\nabla_{h_{k+1}}\nabla_{h_k}\cdots\nn_{h_1} P_{s,t} f(x)   &= \nabla_{h_{k+1}}
\mE\bigg[ (P_{s_{k},t}f)(X_{s,s_k}(x))\prod_{i=1}^{k}\xi_{s_{i-1},s_{i}}^{\breve h_{i}} \bigg]\\
&=\mE\bigg[ \nabla_{h_{k+1}}\big(P_{s_k,t}f(X_{s,s_k}(x))\big)\prod_{i=1}^{k}\xi_{s_{i-1},s_{i}}^{\breve h_i}\bigg]\\
 &= \mE\bigg[ \big(\nabla_{\breve h_{k+1}}P_{s_k,t}f\big)(X_{s,s_k}(x))\prod_{i=1}^{k}\xi_{s_{i-1},s_{i}}^{\breve h_{i}}\bigg]\\
&= \mE\bigg[f(X_{s,t}(x))\prod_{i=1}^{k+1}\xi_{s_{i-1},s_{i}}^{\breve h_{i}}\bigg],
\end{align*}
where in the last step we have used the independence of $\big\{X_{s,s_k}(x), \xi_{s_{i-1},s_{i}}^{\breve h_{i}}, i=1,\cdots,k\big\}$ and
$\big\{X_{s_k, t}(x), \xi_{s_{k},s_{k+1}}^{\breve h_{k+1}}\big\}$.
The proof is complete.
\end{proof}

\bl\label{LNN2}
 For any $p\ge 1$ and $\phi\in\sU$, there is a constant $C=C(\phi, p,\kappa)>0$, where $\kappa$ is given in \eqref{IV}, such that for all $0\leq s<t<\infty$,
\begin{align}\label{AA9}
\big\|\phi\big(\big|X^{(1)}_{s,t}(0)\big|\big)\big\|_p\leq C\phi((t-s)^{\ff 3 2}),\ \
\big\|\phi\big(\big|X^{(2)}_{s,t}(0)\big|\big)\big\|_p\leq C\phi((t-s)^{\ff 1 2}),
\end{align}
where $\|\cdot\|_p:= (\E |\cdot|^p)^{\ff 1 p}.$
\el
\begin{proof}
First of all, by \eqref{SL} and Burkholder's inequality, for any $p\geq 1$  there is a constant $C=C(p,\kappa)>0$ such that for all $0\leq s<t<\infty$,
\begin{align}\label{ER11}
\big\|X^{(1)}_{s,t}(0)\big\|_p\leq C(t-s)^{\ff 3 2},\ \ \big\|X^{(2)}_{s,t}(0)\big\|_p\leq C(t-s)^{\ff 1 2}.
\end{align}
On the other hand, since $\phi\in\sU$ is increasing,   by \eqref{UU0}   we obtain
\begin{align*}
\big\|\phi\big(\big|X^{(1)}_{s,t}(0)\big|\big)\big\|_p &=\big\|\phi\big((t-s)^{\ff 3 2}| (t-s)^{-\ff 3 2}X^{(1)}_{s,t}(0) | \big)\big\|_p\\
&\le C\phi\big((t-s)^{\ff 3 2}\big) \big\|  1+ |(t-s)^{-\ff 3 2}X^{(1)}_{s,t}(0) |^\dd\big\|_p.
\end{align*}
Combining this with \eqref{ER11}
we prove the first estimate. Similarly, we can prove the second estimate.
\end{proof}

Below we present a simple consequence of the above formula, which will play crucial roles in the next section.
In particular, as in \cite{Ch-So-Zh}, the pointwise estimate results given below
allow us to borrow the H\"older regularity of $b^{(1)}$ to compensate the singularity along the first direction induced by the degeneracy.

\bc\label{C2.7} Let   $\phi,\psi\in\sU$.  For any $T>0$ and $m,k\in \mN_0=:\{0\}\cup\mN$,
there exists a constant $C >0$ such that for any $0\leq s<t\leq T$ and any constants $K_1,K_2\ge 0$,
\beq\label{Es22}\beg{split}
\|(\nabla^{(1)})^{\otimes m} (\nabla^{(2)})^{\otimes k} P_{s,t}f\|(0) \leq C \big(K_1\phi((t-s)^{\frac{3}{2}})+K_2\psi((t-s)^{\frac{1}{2}})\big)
(t-s)^{-\frac{3m}2-\ff k 2}
\end{split}\end{equation}
holds for any measurable function $f$ on $\R^{d_1+d_2} $ satisfying
\begin{align}\label{FF}
|f(x)|\le K_1 \phi(|x^{(1)}|)+K_2 \psi(|x^{(2)}|).
\end{align}
Consequently, for any $m\in\mN$, $k\in\mN_0$ and any measurable function $f$ on $\R^{d_1+d_2}$,
\begin{align}
&\quad\|(\nabla^{(1)})^{\otimes m}(\nabla^{(2)})^{\otimes k} P_{s,t} f\|_\infty
\leq C[f]_{\phi,\infty}\phi((t-s)^{\frac{3}{2}})(t-s)^{-\frac{3m}2-\ff k 2},\label{Es2}\\
&\|(\nabla^{(2)})^{\otimes k} P_{s,t} f\|_\infty
\leq C\big([f]_{\phi,\infty}\phi((t-s)^{\frac{3}{2}})+[f]_{\infty,\psi}\psi((t-s)^{\ff 1 2})\big)(t-s)^{-\ff k 2}.\label{Es02}
 \end{align}
\ec
\begin{proof} We introduce the following notations:
$$
\breve\xi_{s,t}^h=\xi_{s,t}^{\breve h},\  h\in\mR^{d_1+d_2},\ \breve\xi^{(\cdot, 0)}_{s,t}=\big(\breve\xi_{s,t}^{(e_i,0)}\big)_{i=1,\cdots,d_1},
$$
where $\xi^{\breve h}_{s,t}$ is defined by \eqref{xi}, and $(e_i)_{i=1,\cdots, d_1}$ is the standard basis of $\mR^{d_1}$.
Similarly, we can define $\breve\xi^{(0,\cdot)}_{s,t}\in\mR^{d_2}$.
By \eqref{EST}, \eqref{xi} and Burkholder's inequality, we have for any $T>0$ and $p\ge 1$,
\begin{align}\label{ER1}
\|\breve\xi_{s,t}^{(\cdot, 0)} \|_p \ppp (t-s)^{-\ff 3 2},\  \|\breve\xi_{s,t}^{(0,\cdot)} \|_p \ppp (t-s)^{-\ff 1 2},\  0\le s<t\le T,
\end{align}
where $\|\cdot\|_p:= (\E |\cdot|^p)^{\ff 1 p}.$

Let $s_i=s+(t-s)i/(m+k), i=0,1,\cdots, m+k$ be the uniform partition of $[s,t]$.
Using the above notations, by \eqref{BS} we have
$$
\|(\nabla^{(1)})^{\otimes m}(\nabla^{(2)})^{\otimes k} P_{s,t} f\|(0)\le
\E\left\{\big|f(X_{s,t}(0))\big|\cdot\left\|\prod_{i=1}^{m}\breve\xi_{s_{i-1},s_{i}}^{(\cdot,0)}\cdot\prod_{j=m+1}^{m+k}\breve\xi_{s_{j-1},s_j}^{(0,\cdot)}\right\|\right\}.
$$
Estimate \eqref{Es22} follows by H\"older's inequality and \eqref{FF}, \eqref{AA9}, \eqref{ER1}.

In general, for fixed $x_0\in \R^{d_1+d_2}$, let
\begin{align*}
g_{x_0}(x)&:=f\big(x_0^{(1)}+\Gamma_{s,t}x^{(2)}_0, x^{(2)}+x_0^{(2)}\big),\\
f_{x_0}(x)&:=f\big(x^{(1)}+x^{(1)}_0+\Gamma_{s,t}x^{(2)}_0, \ x^{(2)}+x_0^{(2)}\big)-g_{x_0}(x).
\end{align*}
Noticing that $\nn^{(1)}P_{s,t} g_{x_0}\equiv0$, we have
$$(\nabla^{(1)})^{\otimes m}(\nabla^{(2)})^{\otimes k} P_{s,t} f(x_0) = (\nabla^{(1)})^{\otimes m}(\nabla^{(2)})^{\otimes k}  P_{s,t}f_{x_0}(0),\ \ m\not=0.$$
Thus,  \eqref{Es2} follows from \eqref{Es22} with $K_2=0$. As for \eqref{Es02}, it follows by \eqref{Es22}.
\end{proof}

\section{A study for degenerate parabolic equations}

Throughout this section, we fix $T, \lambda>0$ and consider the following degenerate parabolic equation with H\"older coefficients:
\begin{align}
\p_t u_t=\sL^{\Sigma,b}_tu_t-\lambda u_t+f_t,\ u_0=0,\ \ t\in [0,T], \label{Eq1}
\end{align}
where $\sL^{\Sigma,b}_t$ is defined by \eqref{LG} and $f: [0,T]\times \R^{d_1+d_2}\to \R$ is measurable.
 The solution will be used in Section 4 to construct the diffeomorphism on $\R^{d_1+d_2} $ which transforms the original \eqref{SDE} into an equation with regular enough coefficients so that the existence and uniqueness of solutions are proved.

Before studying equation \eqref{Eq1}, we first estimate the gradients on $P_{s,t}(H\cdot\nn^{(i)}f), i=1,2$,
which are nontrivial consequences of Corollary \ref{C2.7}, and will play a crucial role in estimating derivatives of $u_t$ in terms of the formula
\eqref{Es9} below. For fixed $\phi\in\sD_0\cap\sS_0$, let
\begin{align}\label{Lam}
\Lambda^\phi_\lambda(t) =\e^{-\lambda t}t^{-1}\phi(t^{\ff 1 2}),\ \ t\in(0,T].
\end{align}

\subsection{Gradient estimates on $P_{s,t}(H\cdot\nn^{(i)}f)$  }

Below all the constants appearing in $\preceq$ only depends on $T, d_1,d_2$ and $\phi$.

\bl\label{Co213} Let $f\in C^1(\R^{d_1+d_2})$ and $H\in C^1(\R^{d_1+d_2};\R^{d_2})$ with $H(0)=0$. For $0\le s<t\le T$
and $k=0,1$, recalling the definition of $\phi_{[\alpha]}$ in \eqref{DD1}, we have
\begin{align}
&\|(\nn^{(2)})^{\otimes(k+1)} P_{s,t}(H\cdot\nn^{(2)} f)\|(0) \preceq[H]_{\phi} \|\nn^{(2)} f\|_{\infty}\Lambda^\phi_0(t-s),\label{DG2}\\
&\|\nn^{(1)}(\nn^{(2)})^{\otimes k} P_{s,t}(H\cdot\nn^{(2)} f)\|(0) \preceq\|H\|_{\phi_{[(k+1)/3]},\infty}
\|\nn^{(2)} f\|_{\phi_{[(k+1)/3]},\infty}\Lambda^\phi_0(t-s),\label{DG202}\\
&\|\nn^{(1)}(\nn^{(2)})^{\otimes k}P_{s,t} (H\cdot \nn^{(2)} f)\|(0)
  \preceq  [H]_{\phi_{[2/3]}}\Big([f]_{1_{[(k+2)/3]},\infty}+\|\nn^{(2)}f\|_\infty\Big) \Lambda^\phi_0(t-s),\label{DG4}\\
&|\nn^{(1)} P_{s,t} (H\cdot \nn^{(2)} f)|(0) \preceq \big([H]_{1_{[2/3]}}
+\|\nn^{(2)} H\|_{\phi_{[1/9]},\infty}\big) [f]_{1_{[2/3]}}\Lambda^\phi_0(t-s).\label{DG3'}
\end{align}
Moreover, if for some $K>0$,
$$
|H(0,x^{(2)})|\leq K|x^{(2)}|\phi(|x^{(2)}|),
$$
then
\begin{align}
&\|(\nn^{(2)})^{\otimes (k+1)}P_{s,t}(H\cdot\nn^{(1)} f)\|(0) \preceq([H]_{\phi_{[1/3]},\infty}+K)\|\nn^{(1)}f\|_\infty(t-s)^{-\ff k2},\label{DG22}\\
&\|\nn^{(1)}(\nn^{(2)})^{\otimes k}P_{s,t} (H\cdot \nn^{(1)} f)\|(0) \preceq ([H]_{\phi_{[(k+1)/3]},\phi}+K)\|\nn^{(1)} f\|_{1_{[k/3]},\infty}\Lambda^\phi_0(t-s). \label{DG44}
\end{align}
\el

\begin{proof}
{\bf (1)}
Since $H(0)=0$, recalling definition \eqref{DD1} and \eqref{DB88},
we have
$$
|H(x)|\le [H]_{\phi}\phi_{[0]}(|x|),\ \ \|H\cdot\nn^{(2)} f\|_{\phi_{[(k+1)/3]}}\leq\|H\|_{\phi_{[(k+1)/3]}}\|\nn^{(2)} f\|_{\phi_{[(k+1)/3]}}.
$$
So, \eqref{DG2} follows by \eqref{Es22}, and \eqref{DG202} follows by \eqref{Es2}.
\\
\\
{\bf (2)} To prove \eqref{DG4}, we introduce
$$
f_{1}(x):=f(x^{(1)}, 0), \ f_{2}(x):=f(0,x^{(2)}), \ \widehat f_{i}(x):=f(x)-f_{i}(x),\ \ i=1,2.
$$
Moreover, for $\bP_\theta^{(2)}$ being the heat semigroup on $\R^{d_2}$, let
$$H_2^\theta = \bP_\theta^{(2)}H_2 - \bP_\theta^{(2)}H_2(0),\ \ \widehat H_2^\theta= H_2- H_2^\theta,\ \ \theta>0.$$
We have
\beq\label{HY0} H\cdot \nn^{(2)}f =\widehat H_2\cdot \nn^{(2)}f + \widehat H_2^\theta \cdot \nn^{(2)} f+H_2^\theta\cdot \nn^{(2)} f,\ \ \theta\in (0,1].\end{equation}
Below we investigate    these three terms respectively.
\\
\\
{\bf (2a)} Observing from \eqref{DB88} that
$$
|\widehat H_2\cdot\nn^{(2)}f|(x) \le [H]_{\phi_{[(k+1)/3]},\infty}\|\nn^{(2)}f\|_\infty \phi_{[(k+1)/3]}(|x^{(1)}|),$$
by \eqref{Es22} we obtain
\beq\label{HY1}\beg{split}
\|\nn^{(1)}(\nn^{(2)})^{\otimes k}P_{s,t}(\widehat H_2\cdot \nn^{(2)}f)\|(0)\ppp[H]_{\phi_{[(k+1)/3]},\infty} \|\nn^{(2)}f\|_\infty \LL_0^\phi(t-s).
\end{split}\end{equation}
\\
{\bf (2b)} Since by \eqref{AA5} we have $\nn^{(1)}P_{s,g} g=0$ for $g$ depending only on $x^{(2)}$, it follows that
$$\nn^{(1)}P_{s,t} (H_2^\theta\cdot\nn^{(2)}f)= \nn^{(1)}P_{s,t} (H_2^\theta\cdot\nn^{(2)}\widehat f_2)= \nn^{(1)}P_{s,t} \div^{(2)}(\widehat f_2 H_2^\theta) - \nn^{(1)}P_{s,t} (\widehat f_2 \div^{(2)}H_2^\theta).$$
Noting that  Lemma \ref{Le24} and \eqref{DD1}  imply
\beq\label{LST}\beg{split} &|\widehat f_2 H_2^\theta|(x)\le \|\nn ^{(2)}H_2^\theta\|_\infty[f]_{1_{[2/3]},\infty} |x^{(2)}|(|x^{(1)}|+|x^{(1)}|^{\ff 2 3}),\\
&|\widehat f_2 \div^{(2)} H_2^\theta|(x)\le \|\nn ^{(2)}H_2^\theta\|_\infty [f]_{1_{[2/3]},\infty} (|x^{(1)}|+|x^{(1)}|^{\ff 2 3}),
\end{split}\end{equation}
from  \eqref{AA5} and Corollary \ref{C2.7}  we obtain
\beq\label{LST2}\beg{split} &\|\nn^{(1)} P_{s,t}(H_2^\theta \cdot \nn^{(2)}f)\|(0)\le \|\nn^{(1)}  P_{s,t}\div^{(2)}(\widehat f_2 H_2^\theta )\|(0)
+ \|\nn^{(1)}  P_{s,t}(\widehat f_2 \div^{(2)}H_2^\theta) \|(0)\\
&\qquad\qquad\le \|\nn^{(1)} \nn^{(2)} P_{s,t} (\widehat f_2 H_2^\theta )\|(0)+\|\GG_{s,t}\|\cdot\|\nn^{(1)} \nn^{(1)} P_{s,t} (\widehat f_2 H_2^\theta )\|(0)\\
&\qquad\qquad+ \|\nn^{(1)}  P_{s,t}(\widehat f_2 \div^{(2)}H_2^\theta) \|(0)
 \ppp [f]_{1_{[2/3]},\infty} \|\nn^{(2)}H_2^\theta\|_\infty  (t-s)^{-\ff 1 2}.\end{split} \end{equation}
Similarly, using
\beg{align*} & |\widehat f_2 H_2^\theta|(x)\le \|\nn ^{(2)}H_2^\theta\|_\infty [f]_{1_{[1]},\infty} |x^{(2)}|\cdot |x^{(1)}| ,\\
&|\widehat f_2 \div^{(2)} H_2^\theta|(x)\le \|\nn ^{(2)}H_2^\theta\|_\infty [f]_{1_{[1]},\infty} |x^{(1)}|,
\end{align*}
to replace \eqref{LST}, we have
\beq\label{LST3}\|\nn^{(1)} \nn^{(2)} P_{s,t}(H_2^\theta \cdot \nn^{(2)}f)\|(0)\ppp [f]_{1_{[1]},\infty} \|\nn^{(2)}H_2^\theta\|_\infty (t-s)^{-\ff 1 2}.\end{equation}
Moreover, by \eqref{Ed1'},
$$
\|\nn^{(2)}H_2^\theta\|_\infty \ppp [H_2]_{\phi_{[2/3]}}\theta^{-\ff 1 6}\phi(\theta^{\ff 1 2}),\ \ \theta\in (0,1].
$$
Then \eqref{LST2} and \eqref{LST3} yield
\beg{align}
\|\nn^{(1)}(\nn^{(2)})^{\otimes k} P_{s,t}(H_2^\theta \cdot \nn^{(2)}f)\|(0)&\ppp  [f]_{1_{[(k+2)/3]},\infty} [H]_{\infty,\phi_{[2/3]}}
\theta^{-\ff 1 6}\phi(\theta^{\ff 1 2})(t-s)^{-\ff 1 2}. \label{HY2}\end{align}
\\
{\bf (2c)} Since $H(0)=0$, by \eqref{Ed1'} we have
\begin{align*}
|\widehat H_2^\theta(x)|& =\left|\int_0^\theta \big(\pp_r \bP_r^{(2)} H_2(0) - \pp_r \bP_r^{(2)} H_2(x)\big)\d r\right|\\
&\leq 2\int_0^\theta\|\pp_r \bP_r^{(2)} H_2\|_\infty\dif r\ppp [H]_{\infty, \phi_{[2/3]}}\int_0^\theta r^{-\ff 2 3}\phi(r^{\ff 1 2})\dif r\\
&\ppp[H]_{\infty, \phi_{[2/3]}} \theta^{\ff 1 3}\phi(\theta^{\ff 1 2}),\ \ \theta\in(0,1].
\end{align*}
Thus,  it follows from  Corollary \ref{C2.7} that for $\theta\in (0,1]$,
\beg{align}
 \|\nn^{(1)}(\nn^{(2)})^{\otimes k} P_{s,t}(\widehat H_2^\theta \cdot \nn^{(2)}f)\|(0)
 \ppp \|\nn^{(2)}f\|_\infty  [H]_{\infty, \phi_{[2/3]}}
 (t-s)^{-\frac{3+k}{2}}\theta^{\ff 13}\phi(\theta^{\ff 1 2}). \label{HY3}\end{align}
Taking $\theta=(t-s)^3$,   by combining \eqref{HY0} with \eqref{HY1}, \eqref{HY2}  and \eqref{HY3}, we prove   \eqref{DG4}.
\\
\\
 {\bf (3)} We now prove \eqref{DG3'}. Since
  $\nn^{(2)} f_1=0$ and $\nn^{(1)} P_{s,t} (H_2\cdot \nn^{(2)}f_2)=0$, we have
   \beq\label{EN1}\beg{split}
   &\nabla^{(1)}P_{s,t}\big(H\cdot\nabla^{(2)}f\big)=\nn^{(1)}P_{s,t} (H_2\cdot \nn^{(2)}\widehat f_2) + \nn^{(1)}P_{s,t} (\widehat H_2\cdot \nn^{(2)}\widehat f_1)\\
 &=\nabla^{(1)}P_{s,t}\div^{(2)}\big(\widehat f_2H_2+ \widehat f_1\widehat H_2\big)  -\nabla^{(1)}P_{s,t} (\widehat f_2\div^{(2)} H_2 +\widehat f_1 \div^{(2)}\widehat H_2).
  \end{split}\end{equation} Below we estimate these two terms respectively.

Firstly, by  $\div^{(2)} \widehat H_2 (x)= \div^{(2)}H(x) -\div^{(2)}H(0,x^{(2)})$, we have
 \beg{align*} &|\widehat f_2\div^{(2)} H_2|(x)
  \ppp     [f]_{1_{[\ff 2 3]},\infty}\|\nn^{(2)} H\|_\infty(|x^{(1)}|^{\ff 2 3}+|x^{(1)}|),\\
&|\widehat f_1 \div^{(2)}\widehat H_2|(x) \ppp [f]_{\infty,1_{[\ff 2 3]}} [\nn^{(2)}H]_{\phi_{[\ff 1 9]},\infty} (|x^{(2)}|^{\ff 2 3}+|x^{(2)}|)\phi_{[\ff 1 9]}(|x^{(1)}|).
\end{align*}
So, Corollary \ref{C2.7} implies
\beq\label{EN2} \beg{split}&\big|\nn^{(1)}P_{s,t} \big(\widehat f_2\div^{(2)} H_2\big)\big|(0)\ppp  [f]_{1_{[2/3]},\infty}\|\nn^{(2)} H\|_\infty (t-s)^{-\ff 12},\\
&  \big|\nn^{(1)}P_{s,t} \big(\widehat f_1\div^{(2)} \widehat H_2\big)\big|(0)\ppp [f]_{\infty,1_{[2/3]}}[\nn^{(2)} H]_{\phi_{[1/9]},\infty }  \LL_0^\phi(t-s).
\end{split}\end{equation}
Next,   since
\beg{align*} |\widehat f_2 H_2|(x)+ |\widehat f_1\widehat H_2|(x)\le [H]_{1_{[2/3]}} [f]_{1_{[2/3]}}
(|x^{(2)}|+|x^{(2)}|^{\ff 2 3})  (|x^{(1)}|^{\ff 2 3}+|x^{(1)}| ),\end{align*}
 it follows from \eqref{AA5} and  Corollary \ref{C2.7} that
  \beg{equation*} \beg{split}& |\nn^{(1)}  P_{s,t}  \div^{(2)}(\widehat f_2 H_2)|(0)+ |\nn^{(1)} P_{s,t} \div^{(2)}(\widehat f_1\widehat H_2) |(0)\\
  & \ppp \|\nn^{(1)} \nn^{(2)}P_{s,t} (\widehat f_2 H_2+ \widehat f_1\widehat H_2)\|(0)
  +(t-s)\|\nn^{(1)}  \nn^{(1)} P_{s,t} (\widehat f_2 H_2+ \widehat f_1\widehat H_2)\|(0)\\
 &\ppp [H]_{1_{[2/3]}}[f]_{1_{[2/3]}}(t-s)^{-\ff 2 3}\ppp [H]_{1_{[2/3]}}[f]_{1_{[2/3]}}\LL_0^\phi(t-s).
  \end{split} \end{equation*}
  Combining  this with \eqref{EN1} and  \eqref{EN2}, we prove \eqref{DG3'}.
  \\
  \\
{\bf (4)} Noticing that
$$
|H\cdot \nn^{(1)} f|(x)\leq\Big([H]_{\phi_{[1/3]},\infty}\phi_{[1/3]}(|x^{(1)}|)+K|x^{(2)}|\phi(|x^{(2)}|)\Big)\|\nn^{(1)} f\|_\infty,
$$
by Corollary \ref{C2.7}, we obtain \eqref{DG22}. Let $g:=H\cdot \nn^{(1)} f$. Observing that
\begin{align*}
|\widehat g_2|(x)&\leq K[\nn^{(1)}f]_{1_{[k/3]},\infty}|x^{(2)}|\phi(|x^{(2)}|)(|x^{(1)}|+|x^{(1)}|^{\ff k 3})\\
&\quad+[H]_{\phi_{[(k+1)/3]},\infty}\|\nn^{(1)}f\|_\infty\phi_{[(k+1)/3]}(|x^{(1)}|),
\end{align*}
and $\nabla^{(1)}g_2\equiv0$, by Corollary \ref{C2.7} again, we have
\begin{align*}
&\|\nn^{(1)}(\nn^{(2)})^{\otimes k}P_{s,t} (H\cdot \nn^{(1)} f)\|(0)=\|\nn^{(1)}(\nn^{(2)})^{\otimes k}P_{s,t} \widehat g_2\|(0)\\
&\qquad\preceq ([H]_{\phi_{[(k+1)/3]},\phi}+K)\|\nn^{(1)} f\|_{1_{[k/3]},\infty}\Lambda^\phi_0(t-s).
\end{align*}
The proof is complete.
\end{proof}

\subsection{Smooth solutions and apriori  estimates}
In this subsection, we study the key apriori estimates for the smooth solutions of equation \eqref{Eq1}. To this aim  we assume that
%\ \ \
\begin{align}
\sup_{t\in [0,T]} \big(\|\nn^{\otimes k} b_t\|_\infty+\|\nn^{\otimes k}f_t\|_\infty+ \|\nn^{\otimes k}\si_t\|_\infty +\|\si_t\|_\infty+\|\si_t^{-1}\|_\infty\big)<\infty,\ \ k\in\mN. \label{Ass}
\end{align}
For fixed $\phi\in\sD_0\cap\sS_0$, we introduce the following quantities for later use:
\beq\label{SQ'}\beg{split}
&\bar\sQ_\phi:= \sup_{t\in[0,T]}\Big\{[b^{(1)}_t]_{\phi_{[2/3]},\infty}+\|\nabla^{(2)}b^{(1)}_t\|_{\infty,\phi}
+\|\big([\nn^{(2)} b^{(1)}_t] [\nn^{(2)} b^{(1)}_t]^*\big)^{-1}\|_\infty \\
&\qquad\qquad\qquad
+\|\sigma^{-1}_t\|_\infty+\|\sigma_t\|_{\phi_{[2/3]}}+[b^{(2)}_t]_{\phi_{[2/3]},\phi} \Big\},\\
\end{split}\end{equation} and
\beq\label{SQ}  \beg{split}
\sQ_\phi:= \bar \sQ_\phi+ \sup_{t\in[0,T]}[b^{(2)}_t]_{\phi_{[2/3]},\phi^{7/2}},\ \ \sQ'_\phi:= \bar \sQ_\phi+\sup_{t\in[0,T]}\|\nn^{(2)}\si_t\|_{\phi_{[1/9]},\infty},
\end{split}\end{equation}
where $\phi_{[\alpha]}$ is defined in \eqref{DD1}.
By \eqref{Ass}, these quantities are all finite.

 \

The main result of this section is the following, which is the key
%in the proofs of the main results.
in the proofs of Theorems \ref{TW}-\ref{T1.2'}.

\beg{thm}\label{T3.1}  Under \eqref{Ass}, $\eqref{Eq1}$ has a unique smooth solution $u$ such that for all $t\in[0,T]$,
\begin{align}
\begin{split}\label{NNY4}
&\|\nabla u_t\|_{1_{[1/3]},\infty}+\|\nabla^{(1)}\nabla^{(2)} u_t\|_\infty+\|\nabla^{(2)}\nabla^{(2)} u_t\|_{\phi^{3/2}}\\
&\qquad\qquad\qquad\leq C\int^t_0\e^{-\lambda(t-s)}\ff{\phi((t-s)^{\ff 1 2})}{t-s}  [f_s]_{\phi_{[2/3]},\phi^{7/2}}  \dif s,
\end{split}\\
\begin{split}\label{NNY4'}
&\|\nabla u_t\|_{1_{[1/3]},\infty}+\|\nabla\nabla^{(2)} u_t\|_\infty
\leq C'\int^t_0\e^{-\lambda(t-s)}\ff{\phi((t-s)^{\ff 1 2})}{t-s}  [f_s]_{\phi_{[2/3]},\phi}  \dif s,
\end{split}
\end{align}
where $C=C(\phi,\sQ_\phi)$ and $C'= C'(\phi,\sQ_\phi')$ are  increasing in $\sQ_\phi$ and $\sQ_\phi'$ respectively.
 \end{thm}

\paragraph{Remark 3.1.} We emphasize that the constants  in Theorem \ref{T3.1} are   increasing in $\sQ_\phi$ or $\sQ_\phi'$,  since this property enables us to
   make smooth approximations of relevant functionals in the proof of the main results without changing the constants.

\

We first prove the existence and uniqueness of $u$.

\beg{lem}\label{L3.2}
Assume \eqref{Ass}. Then \eqref{Eq1} has a unique smooth solution $u$ such that
\begin{align}
\sup_{t\in [0,T]}\|\nn^k u_t\|_\infty<\infty,\ \ k\in\mN,\ \ \sup_{(t,x)\in[0,T]\times\mR^{d_1+d_2}}\frac{|u_t(x)|}{1+|x|}\leq C\ll^{-1},\ \ \label{Lin}
\end{align} holds for some constant
  $C$ increasing in $\sup_{t\in[0,T]}\big(\|\ff{|b_t|+|f_t|}{1+|\cdot|}\|_\infty+ \|\si_t\|_\infty\big)$.
\end{lem}
\begin{proof}
Let $X_{t,s}(x)=X_{t,s}$ solve the following SDE:
$$
\dif X_{t,s}=b_{T-s}(X_{t,s})\dif s+(0,\sigma_{T-s}(X_{t,s})\dif W_s), \ \ X_{t,t}=x\in \mR^{d_1+d_2},\ \ s\in[t,T].
$$
Notice that $u_{T-t}(x)$ solves the following backward equation:
$$
\p_t u_{T-t}+\sL^b_{T-t}u_{T-t}-\lambda u_{T-t}+f_{T-t}=0.
$$
It is well-known that $u_{T-t}(x)$ has the following probabilistic representation (for example, see \cite[Theorem 4.4]{Zh2}),
$$
u_{T-t}(x)=\int^T_t\e^{\lambda(t-s)} \mE f_{T-s}(X_{t,s}(x))\dif s.
$$
By \eqref{Ass}, we have
$$\sup_{s\in [0,T]} \Big(\|\nn^k f_{T-s}\|_\infty + \big\|\E \|\nn^k X_{t,s}(\cdot)\|\big\|_\infty\Big)<\infty,\ \ k\ge 1.$$
Then $u_{t}$ has  bounded derivatives uniformly in $t\in [0,T]$.   Moreover, by the linear growth of $b$ and $f$, it is easy to derive
the second inequality in \eqref{Lin}.
\end{proof}

In order to prove \eqref{NNY4} and \eqref{NNY4'}, we need the following three lemmas, which will be proved in the next subsection.
\bl\label{Le35} Assume \eqref{Ass}.
\beg{enumerate}\item[$(1)$] There exists a constant $\bar C=\bar C(\phi,\bar\sQ_\phi)$  increasing in $\bar\sQ_\phi$ such that  for any $0\le s<t\le T$,
\begin{align}
 \label{ZZ1}
&\|\nn^{(2)}u_t\|_\infty+
\|\nn^{(2)} \nn^{(2)}u_t\|_\infty\le\bar C \int_0^t \Lambda^\phi_\lambda(t-s)\Big(\|\nn u_s\|_\infty +[f_s]_{\phi}\Big)\d s,
\end{align}
%\item[$(2)$] There exists a constant $C=C(\phi, \sQ_\phi\land \sQ_\phi')$  increasing in $ \sQ_\phi\land \sQ_\phi'$ such that  for any $0\le s<t\le T$,
and for $k=0,1$,
\beq\label{ZZ0}\beg{split}
& \|\nn^{(1)}(\nn^{(2)})^{\otimes k}u_t\|_\infty\le \bar C \int_0^t \Lambda^\phi_\lambda(t-s)
\Big(\|\nabla^{(1)}u_s\|_{1_{[k/3]}, \infty}\\
&\qquad\qquad\qquad\qquad+\|\nn^{(2)} u_s\|_{1_{[(k+2)/3]},\infty}+[f_s]_{\phi_{[(k+1)/3]},\phi}\Big)\d s.
\end{split}
\end{equation}
\item[$(3)$] There exists a constant $ C'=  C'(\phi,\sQ_\phi')$  increasing in $ \sQ_\phi'$ such that  for any $0\le s<t\le T$,
\beq\label{ZZ2'}
\|\nn^{(1)} u_t\|_\infty\le   C'  \int_0^t \Lambda^\phi_\lambda(t-s)\Big(\|\nn^{(2)} u_s\|_{1_{[2/3]}} +[f_s]_{\phi_{[1/3]},\infty }\Big)\d s,
\end{equation}
 \end{enumerate}
\el

\bl\label{Le36} Assume \eqref{Ass}.
There exist constants $C=C(\phi,\sQ_\phi)$ and $C'= C'(\phi,  \sQ_\phi')$ which are  increasing in $\sQ_\phi$ and $ \sQ_\phi'$ respectively, such that for all $0\le s<t\le T$,
\beg{equation}  \label{SS1}
\beg{split}
 \|\nabla^{(1)} u_t\|_{1_{[\ff 1 3]},\infty}
&\le C  \int^t_0\Lambda^\phi_\lambda(t-s)\Big(\|\nn\nabla^{(2)} u_s\|_{\infty}
+\|\nn^{(2)}\nn^{(2)}u_s\|_{\infty,\phi^{3/2}}+[f_s]_{\phi_{[2/3]},\phi^{2}}\Big)\dif s
\end{split}\end{equation} and
\beg{equation}  \label{SS1'}
 \|\nabla^{(1)} u_t\|_{1_{[\ff 1 3]},\infty}
 \le  C'  \int^t_0\Lambda^\phi_\lambda(t-s)\Big(\|\nn\nabla^{(2)} u_s\|_{\infty}+[f_s]_{\phi_{[2/3]},\phi}\Big)\dif s.
\end{equation}
\el

\bl \label{Le37} Assume \eqref{Ass}.
There exists a constant $C=C(\phi,\sQ_\phi)$  increasing in $\sQ_\phi$ such that  for any $0\le s<t\le T$,
\begin{align}
\|\nn^{(2)} \nn^{(2)}u_t\|_{\phi^{3/2}}\leq C\int_0^t \Lambda^\phi_\lambda(t-s)
\Big(\|\nn u_s\|_{\phi^{5/2}} +[f_s]_{\phi^{7/2}}\Big)\d s.\label{ZZ3}
\end{align}
\el

Now we can give

\beg{proof}[Proof of Theorem $\ref{T3.1}$]
Letting
$$
h(t):=\|\nabla u_t\|_{1_{[1/3]},\infty}+\|\nn^{(1)}\nabla^{(2)} u_t\|_{\infty}+\|\nn^{(2)}\nn^{(2)}u_t\|_{\phi^{3/2}},
$$
and combining \eqref{ZZ1}, \eqref{ZZ0}, \eqref{SS1} and \eqref{ZZ3}, we obtain
\begin{align*}
h(t)&\preceq \int^t_0\Lambda^\phi_\lambda(t-s) \Big(h(s)+[f_s]_{\phi_{[2/3]},\phi^{7/2}}\Big)\dif s\\
&=\int_0^t  \e^{-\ll(t-s)}\ff{\phi((t-s)^{\ff 1 2})}{t-s}\Big(h(s)+[f_s]_{\phi_{[2/3]},\phi^{7/2}}\Big)\dif s,
\end{align*}
which yields \eqref{NNY4} by Lemma \ref{Le23}.

Similarly,  \eqref{NNY4'} follows   by combining
\eqref{ZZ0}, \eqref{ZZ2'} and \eqref{SS1'}.
\end{proof}

%\subsection{Proofs of Lemmas \ref{Le35}--\ref{Le37} using freezing equations and Duhamel's representation}
\subsection{Proofs of Lemmas \ref{Le35}--\ref{Le37} by using freezing equations and Duhamel's representation}

To prove Lemmas \ref{Le35}-\ref{Le37} by using results presented in Section 2, we  need to represent $u$ by using $P_{s,t}$. To this end,
 we  introduce the following scheme of freezing coefficients at a fixed point $x_0=(x^{(1)}_0, x^{(2)}_0)\in\mR^{d_1+d_2}$.

Let $y_t$ be the unique solution of the following ODE:
\beq\label{WWG}\ff{\d y_t}{\d t} =- b_t(y_t),\ \ y_0=x_0\in\R^{d_1+d_2}.
\end{equation}
Since $b$ is smooth and has bounded derivatives due to \eqref{Ass},
\begin{align}
\theta_t: x_0\mapsto  y_t \ \text{is\ a\ diffeomorphism\ on}\   \mR^{d_1+d_2}.\label{Dif}
\end{align}
Let $\sL^{x_0}_t$ be the freezing operator defined by
$$
\sL^{x_0}_t u=\tr\Big (A_t\cdot\nabla^{(2)}\nabla^{(2)} u\Big)+(B_t x^{(2)})\cdot\nabla^{(1)} u,
$$
where $A_t:=\Sigma_t(y_t)$ and $B_t:=(\nabla^{(2)}b^{(1)}_t)(y_t)$. Set
$$
\tilde u_t(x)=u_t(x+y_t),\ \tilde f_t(x)=f_t(x+y_t),\ \tilde\Sigma_t(x)=\Sigma_t(x+y_t)-\Sigma_t(y_t),
$$
and
\begin{align*}
\tilde b^{(2)}_t(x)=b^{(2)}_t(x+y_t)-b^{(2)}_t(y_t),\ \tilde b^{(1)}_t(x)=b^{(1)}_t(x+y_t)-b^{(1)}_t(y_t)-\nabla^{(2)} b^{(1)}_t(y_t)x^{(2)}.
\end{align*}
From \eqref{Eq1} and \eqref{WWG}  it is easy to see that $\tilde u$ satisfies
$$
\p_t\tilde u=\sL^{x_0}_t\tilde u-\lambda\tilde u+\tr\big(\tilde\Sigma_t\cdot\nabla^{(2)}\nabla^{(2)}\tilde u\big)+\tilde b\cdot\nabla\tilde u+\tilde f,\ \tilde u_0=0.
$$
Let $P_{s,t}$ be the semigroup generated by $\sL^{x_0}_t$.
By Duhamel's formula, we have
\begin{align}
\tilde u_t&=
\int^t_0\e^{-\lambda(t-s)}P_{s,t}\big(\tr\big(\tilde\Sigma_s\cdot\nabla^{(2)}\nabla^{(2)}\tilde u_s\big)+\tilde b_s\cdot\nabla\tilde u_s+\tilde f_s\big)\dif s.\label{Es9}
\end{align}

Note from the definition of $\tilde b^{(1)}_t(x)$ that
\begin{equation}\label{Es600}
\begin{split}
&|\tilde b^{(1)}_t(0,x^{(2)})|=\big|b^{(1)}_t\big(y^{(1)}_t,x^{(2)}+y^{(2)}_t\big)-b^{(1)}_t(y_t)-\nabla^{(2)} b^{(1)}_t(y_t)x^{(2)}\big|\\
&\qquad\qquad\quad\leq|x^{(2)}|\int^1_0\big|\nabla^{(2)}b^{(1)}_t\big(y^{(1)}_t, rx^{(2)}+y^{(2)}_t\big)-\nabla^{(2)} b^{(1)}_t(y_t)\big|\dif r\\
&\qquad\qquad\quad\leq C [\nabla^{(2)}b^{(1)}_t]_{\infty,\phi}|x^{(2)}|\phi_{[0]}(|x^{(2)}|).
\end{split}
\end{equation}
Combining this with \eqref{SQ} and \eqref{SQ'}, we are able to apply   \eqref{DG2}, \eqref{DG202}, \eqref{DG22} and \eqref{DG44} to derive the following lemma.

\bl\label{Le32} Assume \eqref{Ass}. There exist constants $\bar C=\bar C(\phi,\bar \sQ_\phi)$ increasing in $\bar\sQ_\phi$, such that
for all $0\le s<t\le T$ and $k=0,1$,
\begin{align}
&\big\|\nabla^{(1)}(\nabla^{(2)})^{\otimes k}P_{s,t}\big(\tilde b_s\cdot\nabla\tilde u_s\big)\big\|(0)
\le \bar C\Lambda^\phi_0(t-s)\Big(\|\nabla^{(1)}u_s\|_{1_{[\frac{k}{3}]}, \infty}+\|\nabla^{(2)}u_s\|_{\phi_{[\frac{k+1}{3}]},\infty}\Big),\label{BB3}\\
&\qquad\qquad\big\|(\nabla^{(2)})^{\otimes (k+1)}P_{s,t}\big(\tilde b_s\cdot\nabla\tilde u_s\big)\big\|(0)
\le \bar C \Lambda^\phi_0(t-s)\|\nabla u_s\|_\infty.\label{BB4}
\end{align}
\el

The following lemma is an easy consequence of \eqref{Es2} and \eqref{Es02}.
\bl
There is a constant $C=C(\phi, T)>0$ such that for all $0\le s<t\le T$ and $k=0,1$,
\begin{align}
\|\nabla^{(1)} (\nabla^{(2)})^{\otimes k}P_{s,t}\tilde f_s\|_\infty&\le C \Lambda^\phi_0(t-s)[f_s]_{\phi_{[(k+1)/3]},\infty},\label{WZ2}\\
\|(\nabla^{(2)})^{\otimes (k+1)}P_{s,t}\tilde f_s\|_\infty&\le C \Lambda^\phi_0(t-s)[f_s]_{\phi}.\label{WZ3}
\end{align}
\el

Moreover, by \eqref{DG2}, \eqref{DG4} and \eqref{DG3'}, we have

\bl\label{Le33} Assume \eqref{Ass}.
There exist constants $\bar C=\bar C(\phi,\bar \sQ_\phi)$
and $C'=  C'(\phi,\tt \sQ_\phi)$ which are  increasing in $\bar \sQ_\phi$ and $\sQ_\phi'$ respectively, such that for all $0\le s<t\le T$ and $k=0,1$,
\begin{align}
\begin{split}\label{WZ5}
&\big\|\nabla^{(1)}(\nabla^{(2)})^{\otimes k}P_{s,t}\big(\tr\big(\tilde\Sigma_s\cdot\nabla^{(2)}\nabla^{(2)}\tilde u_s\big)\big)\big\|(0)\\
&\qquad\qquad\le \bar C \Big( \|\nn^{(2)} u_s\|_{1_{[(k+2)/3]},\infty}+\|\nn^{(2)}\nabla^{(2)}u_s\|_\infty\Big)\Lambda^\phi_0(t-s),\\
\end{split}\\
\begin{split}
&\big\|(\nabla^{(2)})^{\otimes (k+1)}P_{s,t}\big(\tr\big(\tilde\Sigma_s\cdot\nabla^{(2)}\nabla^{(2)}\tilde u_s\big)\big)\big\|(0)
 \le \bar C \|\nabla^{(2)}\nn^{(2)}u_s\|_{\infty}\Lambda^\phi_0(t-s),\label{WZ6}
 \end{split}
\end{align}
and
\beq\label{WZ4'}  \big|\nabla^{(1)} P_{s,t}\big(\tr\big(\tilde\Sigma_s\cdot\nabla^{(2)}\nabla^{(2)}\tilde u_s\big)\big)\big|(0)
  \le  C' \|\nn^{(2)} u_s\|_{1_{[2/3]}}\Lambda^\phi_0(t-s).\end{equation}
\el

\

Now we are in a position to give the proofs of Lemmas \ref{Le35}-\ref{Le37}.

\begin{proof}[Proof of Lemma \ref{Le35}] Now, substituting estimates in Lemmas \ref{Le32}-\ref{Le33} into \eqref{Es9}, and
noting that $\tt u_t= u_t(\cdot+y_t)$ where, according to \eqref{Dif},
$y_t$ runs all over $\R^{d_1+d_2}$ as $x_0$ does, and by Lemma \ref{Le23} and  \eqref{Es9},
estimate \eqref{ZZ1} follows from  \eqref{BB4}, \eqref{WZ3} and \eqref{WZ6};
estimate \eqref{ZZ0} follows from  \eqref{BB3}, \eqref{WZ2}, \eqref{WZ5} and  \eqref{ZZ1}; and finally, estimate \eqref{ZZ2'} follows from \eqref{BB3}, \eqref{WZ2} and \eqref{WZ4'}.
\end{proof}

\begin{proof}[Proof of Lemma \ref{Le36}]
For simplicity, constants $C$ and $ C'$ below are corresponding to $\sQ_\phi$ and $ \sQ_\phi'$
respectively as in the statement, which may vary from line to line.
\\
\\
{\bf (1)}
Let $\bP^{(1)}_\theta$ be defined by \eqref{PP11}. Let $w^\theta_t(x):=\p_\theta\bP^{(1)}_\theta u_t(x)$ and
\beg{align*}
g^\theta_t(x)&:=\p_\theta \bP^{(1)}_\theta(b_t\cdot\nabla u_t)(x)-(b_t\cdot\nabla \p_\theta \bP^{(1)}_\theta u_t)(x)+\p_\theta \bP^{(1)}_\theta f_t(x)\\
&+{\rm tr}\big(\pp_\theta \bP^{(1)}_\theta (\Sigma_t \cdot \nn^{(2)}\nn^{(2)} u_t)- \Sigma_t \cdot \pp_\theta \bP^{(1)}_\theta\nn^{(2)}\nn^{(2)} u_t\big)(x).
 \end{align*}
By equation \eqref{Eq1}, we have
$$
\p_t w^\theta_t=\scr L_t^{\Sigma, b} w^\theta_t -\ll w_t^\theta +g^\theta_t.
$$
By \eqref{ZZ0} with $k=0$, we have
\begin{align}\label{PQ1}
\|\nn^{(1)} w^\theta_t\|_\infty\le \bar C  \int_0^t \Lambda^\phi_\lambda(t-s)
\Big([\nn^{(2)}w^\theta_s]_{1_{[2/3]},\infty} +[g^\theta_s]_{\phi_{[1/3]},\phi}\Big)\d s.
\end{align}
By the definition of $w_t^\theta$ and using \eqref{CC1} for $g=1$, $\psi=1_{[2/3]}$ and $\phi=1_{[1/3]}$,  we obtain
\beq\label{DST}
[\nn^{(2)} w_s^\theta]_{1_{[2/ 3]}, \infty}\ppp \|\nn^{(1)} \nn^{(2)} u_s\|_{\infty} \theta^{-\ff 56}.
\end{equation}
Next, by Lemma \ref{Le31} with $\psi=\phi_{[\ff 1 3]}$ and $\phi=1_{[1/3]}$, we obtain
\begin{align}\label{DP3'}
[g^\theta_t]_{\phi_{[\ff 1 3]},\infty}\le C
\Big([b_t]_{\phi_{[\ff 2 3]},\infty}\|\nabla u_t\|_\infty+
[\Sigma_t]_{\phi_{[\ff 2 3]},\infty}\|\nabla^{(2)}\nn^{(2)} u_t\|_\infty+[f_t]_{\phi_{[\ff 2 3]},\infty}\Big)\theta^{-\ff 5 6}.
 \end{align}
Moreover, by Lemma \ref{LN9} for $\psi_1=1_{[\ff 2 3]}$ and $\psi_2=\phi^2$, and using $a\wedge c\leq a^{\ff 1 2}c^{\ff 1 2}$ for $a,c>0$, we obtain
\begin{align*}
[\p_\theta \bP^{(1)}_\theta(b_t\cdot\nabla u_t)-b_t\cdot\nabla \p_\theta \bP^{(1)}_\theta u_t]_{\infty, \phi}\le C
[b_t]_{1_{[\ff 2 3]},\phi^2}\|\nabla u_t\|_{\infty,\phi^2} \theta^{-\ff 5 6}
\end{align*}
and
$$
[\p_\theta \bP^{(1)}_\theta f_t]_{\infty, \phi}\le C [f_t]_{1_{[\ff 2 3]}, \phi^2} \theta^{-\ff 5 6}\le C [f_t]_{\phi_{[\ff 2 3]},\phi^2} \theta^{-\ff 5 6}.
$$
Finally, by Lemma \ref{LN9} for $\psi_1=1_{[1]}$ and $\psi_2=\phi^{\ff 32}$, we obtain
$$
[\pp_\theta \bP_\theta^{(1)} (\Sigma_t \cdot \nn^{(2)}\nn^{(2)} u_t)-
\Sigma_t \cdot \pp_\theta \bP_\theta^{(1)}\nn^{(2)}\nn^{(2)} u_t]_{\infty, \phi}\le C \|\nn^{(2)}\nn^{(2)}u_t\|_{\infty,\phi^{3/2}}\theta^{-\ff 5 6}. $$
Therefore,
$$
[g_t^\theta]_{\infty,\phi}\le C \Big(\|\nabla u_t\|_{\infty,\phi^2}+ \|\nn^{(2)}\nn^{(2)}u_t\|_{\infty,\phi^{ 3/2}}+[f_t]_{1_{[\ff 2 3]}, \phi^2}\Big)\theta^{-\ff 5 6}.
$$
Combining this with \eqref{PQ1}, \eqref{DP3'} and \eqref{DST}, and using \eqref{KL2}, we obtain \eqref{SS1}.
\\
\\
{\bf (2)} We now prove \eqref{SS1'}  in the same way. By \eqref{ZZ2'}   for $(w^\theta,g^\theta)$ in place of $(u,f),$ we have
\beq\label{DP3}
\|\nabla w^\theta_t\|_\infty \le   C' \int^t_0\Lambda^\phi_\lambda(t-s)
\Big(\|\nabla^{(2)} w^\theta_s\|_{1_{[\ff 2 3]}}+[g^\theta_s]_{\phi_{[\ff 1 3]},\infty}\Big)\dif s.
 \end{equation}
Due to   \eqref{DST} and\eqref{DP3'},  we only need to estimate $\|\nn^{(2)} w_s^\theta\|_{\infty,1_{[2/3]}}$. By Lemma \ref{LN9} for $g=1$,
$\psi_1=1_{[1]},\psi_2=1_{[1]}$ and $\phi=1_{[2/3]}$, we have
$$
\|\nn^{(2)} w_s^\theta\|_{\infty,1_{[2/3]}} \ppp \|\nn\nabla^{(2)} u_t\|_\infty \theta^{-\ff 5 6}.
$$
This, together with
\eqref{DP3}, \eqref{DP3'}  and \eqref{DST}, yields
\beg{align*}
&\sup_{\theta\in(0,1)}\|\theta^{\ff 5 6}\p_\theta \bP^{(1)}_\theta\nabla u_t\|_\infty
\le   C' \int^t_0\Lambda^\phi_\lambda(t-s)\Big(\|\nn\nabla^{(2)} u_s\|_{\infty}+\|\nabla u_s\|_\infty+[f_s]_{\phi_{[\ff 2 3]},\infty}\Big)\dif s.
\end{align*}
By   Lemma \ref{LN} for $\phi(s)=s^{\ff 1 3}$, this implies \eqref{SS1'}.
  \end{proof}

\begin{proof}[Proof of Lemma \ref{Le37}]
Let $\bP_\theta$ be the semigroup on $\R^{d_1+d_2}$. Let $w^\theta_t=\p_\theta \bP_\theta u_t $ and
\beg{align*}
g^\theta_t(x)&= \p_\theta \bP_\theta(b_t\cdot\nabla u_t)(x)-(b_t\cdot\nabla \p_\theta \bP_\theta u_t)(x)+\p_\theta \bP_\theta f_t(x)\\
& + {\rm tr}\big(\pp_\theta \bP_\theta(\Sigma_t \cdot \nn^{(2)}\nn^{(2)} u_t)- \Sigma_t \cdot \pp_\theta \bP_\theta\nn^{(2)}\nn^{(2)} u_t\big)(x).
 \end{align*}
By equation \eqref{Eq1} we have
$$
\p_t w^\theta_t=\scr L_t^{\Sigma, b} w^\theta_t -\ll w_t^\theta   +  g^\theta_t.
$$
Thus, by \eqref{ZZ1} we have
\beq\label{PQQ}
\|\nn^{(2)}\nn^{(2)} w_t^\theta\|_\infty
 \ppp \int_0^t  \Lambda^\phi_\lambda(t-s) \Big(\|\nn  w_s^\theta\|_\infty+  [ g_s^\theta]_{\phi} \Big)  \d s.
  \end{equation}
On the other hand, by \eqref{Ed1'},  we have
$$
\|\nn w_t^\theta\|_\infty =\|\pp_\theta \bP_{\theta}\nn u_t\|_\infty\ppp \theta^{-1}\phi^{5/2}(\theta^{\ff 1 2})\|\nn u_t\|_{\phi^{5/2}},
$$
and by \eqref{CC1},
$$
[g^\theta_t]_{\phi}\ppp\theta^{-1}\phi^{5/2}(\theta^{\ff 1 2})\Big([b_t]_{\phi^{7/2}}\|\nabla u_t\|_\infty+[f_t]_{\phi^{7/2}}+[\Sigma_t]_{\phi^{7/2}}\|\nn^{(2)}\nn^{(2)}u_t\|_\infty\Big).
$$
Substituting these two estimates into \eqref{PQQ} and noticing that by (ii) of Proposition \ref{Pro1},
$$
\int^t_0s^{-1}\phi^{5/2}(s)\dif s+t\int^1_t s^{-2}\phi^{5/2}(s)\dif s\preceq \phi^{3/2}(t), \ t\in(0,1],
$$
by \eqref{*WF}, we obtain
\begin{align*}
\|\nn^{(2)}\nn^{(2)}u_t\|_{\phi^{3/2}}
\ppp \int_0^t  \Lambda^\phi_\lambda(t-s)\Big(\|\nn u_s\|_{\phi^{5/2}}+  \|\nn^{(2)}\nn^{(2)}u_s\|_\infty+[f_s]_{\phi^{7/2}} \Big)  \d s,
\end{align*}
which gives the desired estimate by Lemma \ref{Le23}.
\end{proof}

\subsection{Classical solutions of \eqref{Eq1}}

In this subsection we prove the existence and stability of classical solutions to equation \eqref{Eq1}.
\beg{thm}\label{T3.10}
Assume $\sQ_\phi<\infty$. For any $f:[0,T]\times\mR^d\to\mR$ with
$$
\sup_{s\in[0,T]}[f_s]_{\phi_{[2/3]},\phi^{7/2}}<\infty,
$$
there exist a unique classical solution $u$ to $\eqref{Eq1}$ such that for all $t\in[0,T]$,
\begin{align}\label{ED101}
\|\nabla u_t\|_{1_{[1/3]},\infty}+\|\nabla^{(2)}\nabla^{(2)} u_t\|_{\phi^{3/2}}
\leq C\int^t_0\e^{-\lambda(t-s)}\ff{\phi((t-s)^{\ff 1 2})}{t-s}  [f_s]_{\phi_{[2/3]},\phi^{7/2}}  \dif s.
\end{align}
Moreover, let $(b^k,\sigma^k, f^k)_{k\in\mN_\infty}$ be a sequence of functions. Let $\sQ^k_\phi$ be defined as in \eqref{SQ} in terms of $(b^k,\sigma^k)$.
Assume that
$$
\sup_k\Big(\sQ^k_\phi+\sup_{s\in[0,T]}[f^k_s]_{\phi_{[2/3]},\phi^{7/2}}\Big)<\infty,
$$
and for each $t>0, x\in\mR^{d_1+d_2}$,
 $$
\lim_{k\to\infty} \|\sigma^k_t(x)-\sigma^\infty_t(x)\|+|b^k_t(x)-b^\infty_t(x)|+|f^k_t(x)-f^\infty_t(x)|=0.
 $$
Let $u^k_t(x)$ be the unique classical solution of \eqref{Eq1} corresponding to $(b^k,\sigma^k, f^k)$ for each $k\in\mN_\infty$.
Then for each $T,R>0$,
\begin{align}\label{CON}
\lim_{k\to\infty}  \sup_{t\in [0,1], |x|\le R}
 \Big(|u_t^k-u^\infty_t|+ |\nn(u_t^k-u^\infty_t)|+\|\nabla^{(2)}\nn^{(2)}(u_t^k-u^\infty_t)\|\Big)(x) =0.
\end{align}
\end{thm}
\begin{proof}
{\bf (1)} Let $\varrho$ be a non-negative smooth function with compact support in $\mR^d$
having
$$
\int_{\mR^d} \varrho (x) \dif x=1.
$$
For $n\in\mN$, define $\varrho_n (x)=n^d \varrho (nx)$ and
\begin{align}
b^n_t=\varrho_n* b_t,\ \ \si^n_t=\varrho_n*\si_t,\ \ f^n_t:=\varrho_n*f_t.\label{EQ10}
\end{align}
Clearly, $b^n, \si^n$ and $f^n$ satisfy \eqref{Ass}. Let $\sQ^n_{\phi}$ be defined by \eqref{SQ} corresponding to $b^n,\sigma^n$.
It is easy to see that for some $n_0$ large enough and all $n\geq n_0$,
$$
\sQ^n_{\phi}\leq 2\sQ_{\phi}.
$$
Let $u^n$ be the unique smooth solution of the following equation
 \beq\label{ED1}
 \pp_t u_t^n = \scr L_t^{\Sigma^n,b^n} u_t^n -\ll u_t^n +f_t^n,\ \ u_0^n=0,\ t\in [0,T],
 \end{equation}
which enjoys the following uniform estimate:
\begin{align}\label{ED11}
\begin{split}
&\|\nabla u^n_t\|_{1_{[1/3]},\infty}+\|\nabla^{(1)}\nabla^{(2)} u^n_t\|_\infty+\|\nabla^{(2)}\nabla^{(2)} u^n_t\|_{\phi^{3/2}}\\
&\qquad\qquad\qquad\leq C\int^t_0\e^{-\lambda(t-s)}\ff{\phi((t-s)^{\ff 1 2})}{t-s}  [f_s]_{\phi_{[2/3]},\phi^{7/2}}  \dif s,
\end{split}
\end{align}
So, Ascoli-Arzela's theorem implies the existence of $u$ such that, up to a subsequence,
$$
\lim_{n\to\infty}  \sup_{t\in [0,1], |x|\le R}
 \Big(|u_t^n-u_t|+ |\nn(u_t^n-u_t)|+\|\nabla^{(2)}\nn^{(2)}(u_t^n-u_t)\|\Big)(x) =0,\ \ R>0.
$$
 By taking limits for \eqref{ED1} and inequality \eqref{ED11}, we obtain the existence of classical solutions of \eqref{Eq1} as well as the estimate \eqref{ED101}.
 \\
 \\
{\bf (2)} We use a contradiction argument. Suppose that \eqref{CON} does not hold. Then there is a subsequence $k_m$ such that
$$
\varliminf_{m\to\infty}  \sup_{t\in [0,1], |x|\le R}
 \Big(|u_t^{k_m}-u^\infty_t|+ |\nn(u_t^{k_m}-u^\infty_t)|+\|\nabla^{(2)}\nn^{(2)}(u_t^{k_m}-u^\infty_t)\|\Big)(x)>0.
$$
On the other hand, repeating the proof in step {\bf (1)}, since $u^\infty$ is the unique solution of \eqref{Eq1} corresponding to $(b^\infty,\sigma^\infty, f^\infty)$, there is a subsubsequence $k'_m$ such that
$$
\lim_{m\to\infty}  \sup_{t\in [0,1], |x|\le R}
 \Big(|u_t^{k'_m}-u^\infty_t|+ |\nn(u_t^{k'_m}-u^\infty_t)|+\|\nabla^{(2)}\nn^{(2)}(u_t^{k'_m}-u^\infty_t)\|\Big)(x) =0.
$$
Thus, we obtain a contradiction, and so, \eqref{CON} holds.
\end{proof}

We also have the following existence of H\"older classical solutions under H\"older assumptions.
\beg{thm}\label{T3.11}
Assume for some $\aa\in (\ff 2 3,1),\bb\in (0,\ff 12)$,
\beq\label{SQ2}\beg{split}& \sQ_{\aa,\bb}:=\sup_{t\in [0,T]} \Big\{[b_t^{(1)}]_{1_{[\alpha]},\infty}+ \|\nn^{(2)} b_t^{(1)}\|_{\infty,1_{[\beta]}}
 + \big\|\big([\nn^{(2)} b^{(1)}_t] [\nn^{(2)} b^{(1)}_t]^*\big)^{-1}\big\|_\infty \\
 &\qquad \qquad \qquad +\|b_t^{(2)}\|_{1_{[\alpha]},1_{[\beta]}}+\|\sigma_t\|_\infty+\|\sigma^{-1}_t\|_\infty +\|\sigma_t\|_{1_{[\alpha]}}\Big\}<\infty.
 \end{split}\end{equation}
Then  for any $\eps\in(0,\beta\wedge(\alpha-\frac{2}{3}))$, there exist a unique solution $u$ to $\eqref{Eq1}$ and constants
$\delta\in (0,1)$ depending only on $\aa,\bb$, and $C=C(\sQ_{\alpha,\beta},\eps,\delta)>0$,
which is increasing in $\sQ_{\aa,\bb}$, and such that for all $t\in[0,T]$ and $\lambda\geq 0$,
\begin{align}\label{NNY0}
\|\nabla u_t\|_{1_{[1/3]}}+\|\nabla\nabla^{(2)}u_t\|_{1_{[\eps]}, 1_{[\eps]}} \leq C\int^t_0\e^{-\lambda(t-s)}  (t-s)^{-\delta}  [f_s]_{1_{[\alpha]},1_{[\beta]}}\dif s.
\end{align}
\end{thm}

\begin{proof}
First of all, we assume \eqref{Ass}.
Following the proof of Lemma \ref{Le37}, by Lemma \ref{Le26}, we have for any $\eps\in(0,\beta\wedge(\alpha-\frac{2}{3}))$,
 \begin{align*}
[g^\theta_t]_{1_{[\alpha-\eps]}, 1_{[\beta-\eps]}}\preceq \Big([b_t]_{1_{[\alpha]}, 1_{[\beta]}}\|\nabla u_t\|_\infty+[f_t]_{1_{[\alpha]}, 1_{[\beta]}}
+[\Sigma_t]_{1_{[\alpha-\eps]}, 1_{[\beta-\eps]}}\|\nn^{(2)}\nn^{(2)} u_t\|_\infty\Big)\theta^{\frac{\eps}{2}-1}.
 \end{align*}
Noticing that for any $\phi\in\sS_0$,
$$
\phi_{[2/3]}(s)=s^{3/2}\phi(s)\preceq1_{[\alpha-\eps]}(s)=s^{\alpha-\eps},\ \ \phi^{7/2}(s)\preceq1_{[\beta-\eps]}(s)=s^{\beta-\eps},\ \ s\in[0,1],
$$
by \eqref{NNY4}, we obtain that for some $\delta\in(0,1)$,
 \begin{align*}
& \|\nabla w^\theta_t\|_{1_{[1/3]},\infty}+\|\nabla\nabla^{(2)} w^\theta_t\|_\infty
\preceq \int^t_0\e^{-\lambda(t-s)}(t-s)^{-\delta}  [g^\theta_s]_{1_{[\alpha-\eps]}, 1_{[\beta-\eps]}} \dif s\\
&\preceq \theta^{\frac{\eps}{2}-1}\int^t_0\e^{-\lambda(t-s)}(t-s)^{-\delta}\Big(\|\nabla u_s\|_\infty
+ \|\nn^{(2)}\nn^{(2)} u_s\|_\infty+ [f_s]_{1_{[\alpha]}, 1_{[\beta]}}\Big)  \dif s,
 \end{align*}
 which in turn gives \eqref{NNY0} by Lemma \ref{LN} and \eqref{NNY4}.
 In general, we can follow the same approximation as done in Theorem \ref{T3.10}.
 \end{proof}
 \section{Proofs of Main Results}

 \begin{proof}[Proof of Theorem \ref{TW}] The existence of weak solution is well known, see e.g.
 \cite[Theorem 2.2 and Remark 2.1, Chapter IV]{Ik-Wa} and \cite{St0}.
 So, we only prove the uniqueness.
Let $(\Omega,\sF,\mP; X_t, W_t)$ and $(\Omega',\sF',\mP'; X'_t, W'_t)$ be two weak solutions of SDE \eqref{SDE} with $X_0=X_0'=x\in \R^{d_1+d_2}$. Fix $T>0$ and $f\in C^\infty_b([0,T]\times\mR^{d_1+d_2})$.
For any $n\ge 1$, let $\si^n$ and $b^n$ be in \eqref{EQ10}, and
let $\sQ_n$ and  $\sQ_n'$ be the numbers defined in \eqref{SQ} and \eqref{SQ'} for  $(b^n,\si^n)$ in place of $(b,\si)$.
It is easy to see that for some $n_0$ large enough and all $n\geq n_0$,
$$
\sQ_n\leq 2\sQ,\ \  \sQ_n'\le  2\sQ'.
$$
By Theorem \ref{T3.1} for $(\scr L_{T-t}^{\Sigma^n,b^n}, f_{T-t})$ in place of $(\scr L_{t}^{\Sigma,b}, f_{t})$, for any $\ll\ge 0$   the equation
 \beq\label{ED} \pp_t  u_t^n = \scr L_{T-t}^{\Sigma^n,b^n}  u_t^n -\ll  u_t^n +f_{T-t},\ \  u_0^n=0, t\in [0,T]
 \end{equation}
 has a unique   solution $u^n: [0,T]\times \R^{d_1+d_2}\to \R$ such that
\beq\label{*GB}
 \sup_{t\in [0,T], n\ge 1} \Big(\|\nn u_t^n\|_{1_{[\ff 1 3]},\infty}
 +\|\nn \nn^{(2)}u_t^n\|_\infty\Big)\le \vv(\ll):= C\int_0^T \e^{-\ll(t-s)}\ff{\phi((t-s)^{\ff 1 2})}{t-s} \d s
\end{equation}
 for some constant $C>0$. So, Ascoli-Arzela's theorem implies the existence of
 $$
 u: [0,T]\times \R^{d_1+d_2}\to \R $$
 such that, up to a subsequence,
 \beq\label{ED0}
 \lim_{n\to\infty} \sup_{t\in [0,T], |x|\le R} \Big(|u_t^n-u_t|+ \|\nn^{(2)}(u_t^n-u_t)\|\Big)(x) =0,\ \ R>0,
 \end{equation}
 and, moreover,
 \beq\label{ED2} \sup_{t\in [0,T]}\Big([u_t]_{1_{[1]}} + [\nn^{(2)}u_t]_{1_{[1]}}\Big)\le\vv(\ll).\end{equation}
Now, due to \eqref{SDE} and \eqref{ED} with $\ll=0$,  It\^o's formula for $u_{T-t}^n(x)$ implies
\beg{align*}
0&=u_T^n(x)+\int^T_0\mE\big\{(\p_s+\sL^{\Sigma, b}_s )u_{T-s}^n(X_s)\big\}\dif s\\
&=u_T^n(x)  +\E\int_0^T\Big\{{\rm tr}\big[\big(\Sigma_s-\Sigma_s^n)\nn^{(2)}\nn^{(2)}u_{T-s}^n\big]
+(b_s-b_s^n)\cdot \nn u_{T-s}^n-f_s\Big\}(X_s)\d s.\end{align*}
So, according to \eqref{*GB}, \eqref{ED0} and noting that $\{|b_t-b_t^n|+\|\si_t-\si_t^n\|\}_{n\ge 1}$ is bounded uniformly in $t\in [0,T]$ and converges to $0$ as $n\to\infty$,
by the dominated convergence theorem, letting $n\to\infty$ we obtain
$$u_T (x)=\int^T_0\mE f_s(X_s)\dif s.
$$
By the same reason, we also have
$$
 u_T(x)=\int^T_0\mE' f_s(X'_s)\dif s.
$$
Hence,
$$
\int^T_0\mE f_s(X_s)\dif s=\int^T_0\mE' f_s(X'_s)\dif s,\ \ f\in C^\infty_b([0,T]\times\mR^{d_1+d_2}).
$$
By \cite[Corollary 6.2.4]{St}, this implies  the weak uniqueness.
\end{proof}

\begin{proof}[Proof of Theorem \ref{T1.1}]
If \eqref{1.7} holds, then the non-explosion and estimate \eqref{1.8} follows by \cite[Lemma 2.2]{Zh3}. So, we only prove the existence and uniqueness of local solutions.
\\
\\
{\bf (1)} We first assume   that  {\bf (A)}  holds for some $C_n=C, \phi_n=\phi$ and $\gg_n=\gg$ independent of $n\ge 1.$
Noting that  $\gg\in \C$ implies $\gg(r)\le c r^{-\ff 1 4}$ for some $c>0$ and all $r\in (0,1]$, in this case we have either  $\sQ_\phi<\infty$ or $\sQ_\phi'<\infty$. Due to the existence of the weak solution as explained in the proof of Theorem \ref{T1.1}, by the Yamada-Watanabe principle  \cite{Ya-Wa}, we only need to prove the pathwise uniqueness.

Let $b^n$, $\sigma^n$ be defined as in \eqref{EQ10}. As in the proof of Theorem \ref{TW}, by Theorem \ref{T3.1}
 for $(\scr L_{T-t}^{\Sigma^n,b^n}, b_{T-t}^n)$ in place of $(\scr L_{t}^{\Sigma,b}, f_{t})$, the equation
 \beq\label{ED'} \pp_t \u_t^n = \scr L_{T-t}^{\Sigma^n,b^n} \u_t^n -\ll \u_t^n +b_{T-t}^n,\ \ \u_0^n=0, t\in [0,T]
 \end{equation}
 has a unique   solution $\u^n: [0,T]\times \R^{d_1+d_2}\to \R^{d_1+d_2}$ such that \eqref{*GB}--\eqref{ED2} hold for   $\u^n$ and some
 $\u: [0,T]\times \R^{d_1+d_2}\to \R^{d_1+d_2}$ in place of $u^n$ and $u$.
 Let
 $$
 \Phi_t(x)= x+\u_{T-t}(x),\ \ \ t\in [0,T], x\in \R^{d_1+d_2}.
 $$
 Then for large enough $\ll>0$, $\Phi_t$ is a homeomorphism on $\R^{d_1+d_2}$ such that
 \beq\label{ED3} \sup_{t\in [0,T]} \Big([\Phi_t]_{1_{[1]}} + [\Phi_t^{-1}]_{1_{[1]}}\Big)<\infty;\end{equation}that is, both $\Phi_t$ and $\Phi_t^{-1}$ are Lipschitz continuous uniformly in $t\in [0,T].$

 Now, if $X_t$ solves \eqref{SDE} up to a stopping time $\tau\le T$, then by It\^o's formula and \eqref{ED'}, we have
%\beg{align*} &\d \big\{X_t +\u_{T-t}^n(X_t)\big\}\\
% &= \Big\{\ll \u_{T-t}^n+{\rm tr}\big[(\Sigma_t- \Sigma_t^n)\nn^{(2)}\nn^{(2)}\u_{T-t}^n\big]+(b_t-b_t^n)\cdot\nn \u_t^n +b-b^n\Big\}(X_t)\d t\\
%&\quad+ \big(0,  \si_t \d W_t\big)+ (\nn_{\si_t\d W_t}^{(2)} \u_{T-t}^n)(X_t),\ \ t\in [0,\tau].
%\end{align*}
% Thus, $\P$-a.s. for all $t\in [0,\tau],$
\beg{align*}
 &X_t + \u_{T-t}^n(X_t)- X_0 - \u_{T}(X_0)\\
  &=  \int_0^t \Big\{\ll \u_{T-s}^n+{\rm tr}\big[(\Sigma_s- \Sigma_s^n)\nn^{(2)}\nn^{(2)}\u_{T-s}^n\big]+(b_s-b_s^n)\cdot\nn \u_{T-s}^n +b_s-b^n_s\Big\}(X_s)\d s\\
  &\quad+ \int_0^t \big(0,  \si_s \d W_s\big)+ \int_0^t (\nn_{\si_s\d W_s}^{(2)} \u_{T-s}^n) (X_s),\  \ t\in [0,\tau],\ \mP\text{-a.s.}
  \end{align*}
 So, as explained in the proof of Theorem \ref{TW},   by letting $n\to\infty$, we obtain for $t\in [0,\tau]$,
 $$
 \Phi_t(X_t) =\Phi_0(X_0) +\int_0^t \ll \u_{T-s}(X_s)\d s +   \int_0^t \big(0,  \si_s (X_s)\d W_s\big)+ \int_0^t (\nn_{\si_s\d W_s}^{(2)} \u_{T-s})(X_s).
 $$
Therefore, if $(X_t)_{t\in [0,\tau]}$ solves \eqref{SDE}, then $Y_t :=\Phi_t(X_t)$ solves the following SDE for $t\in [0,\tau]:$
\beq\label{SDE'}
\d Y_t= \ll (\u_{T-t}\circ\Phi_t^{-1})(Y_t)\d t +  \big\{(\nn_{\si_t\d W_t}^{(2)} \Phi_t)\circ\Phi_t^{-1}\big\}(Y_t).
\end{equation}
Since by \eqref{ED2} and \eqref{ED3},
both $\u_{T-t}\circ \Phi_t^{-1}$ and $(\nn^{(2)} \Phi_t)\circ\Phi_t^{-1}$ are Lipschitz continuous uniformly in $t\in [0,T],$
from the condition  \eqref{LB1'} or \eqref{LB2'} on $\si$ we  see  that \eqref{SDE'} has a unique solution up to time $T$ (see   \cite[Theorem 4.1]{SWY}). So, the pathwise uniqueness of \eqref{SDE} holds up to any stopping time less than $T$. By the arbitrary of $T>0$ we conclude that    \eqref{SDE} has a unique solution for all $t\ge 0$.
\\
\\
{\bf (2)} Next, if $\si(x)$ and $b(x)$ do not depend on $x^{(1)}$, then so does   $\u^n(x)$. In this case, if \eqref{LPT} holds with $\phi_n$ and $\gg_n$ uniformly in $n\ge 1$, then   $\bar\sQ_\phi<\infty$ for some $\phi\in \D_0\cap\sS_0$, so that  by \eqref{ZZ1} we may repeat the above argument to prove   the pathwise uniqueness.
\\
\\
{\bf (3)} In general,  by a localization argument as in  \cite[Proof of Theorem 1.1]{Wa-Zh},
we obtain the local existence and uniqueness of SDE \eqref{SDE} up to explosion time $\zeta$.
More precisely, for any $m\ge 1$,
let $\theta_m\in C_0^\infty(\R^{d_1+d_2};\R^{d_1+d_2})$ be such that $\theta_m(x)=x$ for $|x|\le m$.   Define
\beg{align}\label{CUT}
\si_m(x)= \si\circ\theta_m(x),\ \   b_m^{(2)}(x)= b^{(2)}\circ\theta_m(x),\ \
b_m^{(1)}(x)= b^{(1)}(\theta^{(1)}_m(x), x^{(2)}).
\end{align}
Here and below, for simplicity of notation, we shall drop the time variables in $b$ and $\sigma$ since it does not play any role in the proof.
If   {\bf (A)} or \eqref{LPT} holds, then   for any $m\in\N$,   $\si_m$ and $b_m$ satisfy the same  assumption   for some uniform $C, \phi$ and $\gg$.
For fixed $X_0\in \R^{d_1+d_2}$, let $X_t^m$ with $X_0^m=X_0$  be the unique solution to \eqref{SDE}
for $(\si_m,b_m)$ in place of $(\si, b)$.  Since $b_m(x)=b(x), \si_m(x)=\si(x)$ for $|x|\le m$,   $X_t^m$ solves the original equation \eqref{SDE} up to the stopping time
$$
\tau_m:= \inf\{t\ge 0: |X_t^m|\ge m\}.
$$
By step (1), we have $X_t^n=X_t^m$ for $t\le \tau_n\land\tau_m$, and $\tau_n$ is increasing in $n$. Letting $\zeta=\lim_{n\to\infty}\tau_n$, we see that
$$
X_t:= \sum_{t\in [\tau_{n-1}, \tau_n)}X_t^n,\ \ \tau_0:=0,\ \ t<\zeta
$$
is the unique solution to \eqref{SDE} with life time $\zeta$, i.e. $\limsup_{t\to\zeta}|X_t|=\infty$ holds a.s. on $\{\zeta<\infty\}$.
\end{proof}

\begin{proof}[Proof of Theorem \ref{T1.5}]
{\bf (1)} First of all, we assume that the global conditions in the theorem hold for $(b^k,\sigma^k)_{k\in\mN}$.
In this case, let $\u^k$ be the unique classical solution of \eqref{Eq1} in Theorem \ref{T3.10} corresponding to $(b^k,\sigma^k, b^k)$. Define
$$
 \Phi^k_t(x)= x+\u^k_{T-t}(x),\ \ \ t\in [0,T], x\in \R^{d_1+d_2}.
 $$
 As in the proof of Theorem  \ref{T1.1}, for large enough $\ll>0$, and for each $k\in\mN_\infty$, $\Phi^k_t$ is a homeomorphism on $\R^{d_1+d_2}$ such that
 $$
 \sup_{k\in\mN_\infty}\sup_{t\in [0,T]} \Big([\Phi^k_t]_{1_{[1]}} + [(\Phi_t^k)^{-1}]_{1_{[1]}}\Big)<\infty;
 $$
 By It\^o's formula, $Y^k_t :=\Phi^k_t(X^k_t)$ solves the following SDE for $t\in[0,T]$,
$$
\d Y^k_t= g^k_t(Y^k_t)\d t +  \Theta^k_t(Y^k_t)\d W_t,\ \ Y^k_0=\Phi^k_0(x),
$$
where
$$
g^k_t:=\lambda \u^k_{T-t}\circ(\Phi^{k})^{-1}_t,\ \ \Theta^k_t:= \big(\nn^{(2)}_{\si^k_t\cdot}\Phi^k_t)\circ(\Phi_t^{k})^{-1}.
$$
Moreover, by \eqref{CON}, it is easy to see that for each $t,x\in\mR^d$,
$$
\lim_{k\to\infty}\Big(|g^k_t(x)-g^\infty_t(x)|+\|\Theta^k_t(x)-\Theta^\infty_t(x)\|+|\Phi^k_t(x)-\Phi^\infty_t(x)|\Big)=0.
$$
and by \eqref{ED101}, for all $x,y\in\mR^{d_1+d_2}$,
$$
\sup_{k\in\mN_\infty}\sup_{t\in[0,T]}\Big(|g^k_t(x)-g^k_t(y)|+\|\Theta^k_t(x)-\Theta^\infty_t(y)\|\Big)\leq C|x-y|.
$$
Hence, by \cite[Theorem 15, p.271]{Pro}, we have for each $T,\eps>0$,
$$
\lim_{k\to\infty}\mP\bigg(\sup_{t\in[0,T]}|Y^k_t-Y^\infty_t|\geq\eps\bigg)=0,
$$
which in turn implies \eqref{STA}.
\\
\\
{\bf (2)} In general, by the assumption and \eqref{1.8}, we have the following uniform estimate:
\begin{align}\label{1.88}
\sup_k \E\exp\bigg[\sup_{t\in [0,T]} H(X^k_t(x))^{\vv'} \bigg]
\le \Psi(T) \exp\big[H(x)^\vv\big],\ \ T>0, x\in \R^{d_1+d_2}.
\end{align}
For each $m\in\mN$,
let $\theta_m\in C_0^\infty(\R^{d_1+d_2};\R^{d_1+d_2})$ be such that $\theta_m(x)=x$ for $H(x)\le m$.
Let $\sigma^k_m$ and $b^k_m$ be defined as in \eqref{CUT}, and let $X^{k,m}_t(x)$ be the solution of SDE \eqref{SDE} corresponding to $(\sigma^k_m,b^k_m)$.
Define the stopping times
$$
\tau^k_m:=\inf\{t>0: H(X^k_t(x))\wedge H(X^\infty_t(x))\geq m\}.
$$
Then by \eqref{1.88}, we have
\begin{align}\label{STA1}
\sup_k\mP(\tau^k_m<T)\leq\sup_k\mE \Big(\sup_{t\in[0,T]}H(X^k_t(x))\wedge H(X^\infty_t(x))\Big)/m\to 0,\ \ m\to\infty.
\end{align}
On the other hand, we have
\begin{align*}
\mP\bigg(\sup_{t\in[0,T]}|X^k_t-X^\infty_t|\geq\eps\bigg)&\leq \mP\bigg(\sup_{t\in[0,T]}|X^k_t-X^\infty_t|\geq\eps; \tau^k_m\geq T\bigg)+\mP(\tau^k_m<T)\\
&\leq \mP\bigg(\sup_{t\in[0,T]}|X^{k,m}_t-X^{\infty,m}_t|\geq\eps\bigg)+\mP(\tau^k_m<T),
\end{align*}
which together with step (1) and \eqref{STA1} gives the desired estimate \eqref{STA}.
\\
\\
{\bf (3) }Let $\varphi:\mR_+\to\mR_+$ be a bounded smooth function with $\varphi(r)=r$ for $|r|\leq 1$. Let $\xi^k_t(x):=|X^k_t(x)-X^\infty_t(x)|^2$.
For fixed $R>0$, let $\chi_R:\mR^d\to[0,1]$ be a smooth function with $\chi_R(x)=1$ for $|x|\leq R$ and $\chi_R(x)=0$ for $|x|\geq 2R$.
By Gagliado-Nirenberg's inequality and \eqref{STA3} for some $p>d$, we have
\begin{align*}
&\mE\Big[\sup_{t\in[0,T]}\|\varphi(\xi^k_t)\chi_R\|_\infty\Big]\\
&\leq C\mE\bigg[\sup_{t\in[0,T]}\|\varphi(\xi^k_t)\chi_R\|^{1-\ff d p}_{L^p}\big(\|\varphi(\xi^k_t)\chi_R\|_{L^p}
+\|\chi_R\nn(\varphi(\xi^k_t))\|_{L^p}+1\big)^{\ff d p}\bigg]\\
&\leq C\bigg\{\mE\Big[\sup_{t\in[0,T]}\|\varphi(\xi^k_t)\chi_R\|^{\frac{p(p-d)}{p^2-d}}_p\Big]\bigg\}^{1-d/p^2}\to 0, \ n\to\infty,
\end{align*}
due to \eqref{STA} and the dominated convergence theorem.  So, \eqref{STA4} holds.
\end{proof}

\begin{proof}[Proof of Theorem \ref{T1.2}]
 By Theorem \ref{T1.1}, for each $x\in\mR^{d_1+d_2}$, there is a unique global solution $\{X_t(x),t\geq 0\}$
for SDE \eqref{SDE}. Let $\u_t(x)$ be the unique solution of equation \eqref{Eq1} in Theorem \ref{T3.11}
corresponding to $(\sigma,b)$ and $f=b$.
By \eqref{NNY0}, we have
\begin{align}
\|\nabla \u_t\|_{1_{[1/3]}}+\|\nabla\nabla^{(2)}\u_t\|_{1_{[\eps]}} \leq C\int^t_0\e^{-\lambda(t-s)}  (t-s)^{-\delta}\dif s.\label{JG1}
\end{align}
As in the proof of Theorem \ref{T1.1}, let
$$
\Phi_t(x)=x+\u_{T-t}(x).
$$
By \eqref{JG1}, for large enough $\ll>0$, $\Phi_t$ is a diffeomorphism on $\R^{d_1+d_2}$ such that
 \beq\label{ED13}
 \sup_{t\in [0,T]} \Big(\|\nabla\Phi_t\|_{1_{[\eps]}} +\|\nabla\Phi_t^{-1}\|_{1_{[\eps]}}\Big)<\infty;
 \end{equation}
Moreover, as shown in the proof of Theorem \ref{T1.1} that if $X_t(x)$ solves SDE \eqref{SDE} then $Y_t=\Phi_t(X_t)$ solves \eqref{SDE'}.
By \eqref{ED13} and the condition on $\si$ in Theorem \ref{T1.2}, we have
\begin{align}
\sup_{t\in[0,T]}(\|\nabla (\u_{T-t}\circ\Phi_t^{-1})\|_{1_{[\eps]}}+\|\nn\{(\nn^{(2)}_{\si_t\cdot}\Phi_t)\circ\Phi_t^{-1}\} \|_{1_{[\eps]}})<\infty,\label{EQ22}
\end{align} for some $\vv>0$. So,
by \cite[Theorem 4.6.5]{Ku}, $\{Y_t(\cdot)\}_{t\in[0,T]}$ forms a $C^1$-stochastic diffeomorphism flow, and so does $\{X_t(\cdot):=\Phi_t^{-1}(Y_t(\cdot))\}_{t\in[0,T]}$.
Finally, it is easy to prove \eqref{FL}   from \eqref{SDE'}, \eqref{ED13} and $\sup_{t\in [0,T]}\|\nn\si_t\|_\infty<\infty$.
%below we will prove   the more general result \eqref{FL'}.
\end{proof}

\begin{proof}[Proof of Theorem \ref{T1.2'}]
%{\bf (1)} Let us first prove the apriori estimate  \eqref{FL'}.
As shown in the proof of Theorem \ref{T1.1} that SDE \eqref{SDE1} admits a unique global strong solution $X_t(x)$ and  $Y_t:=\Phi_t(X_t)$ solves (see \eqref{SDE'})
\begin{align}
\dif Y_t=g^a_t(Y_t)\dif t+\Theta_t(Y_t)\dif W_t,\ \ Y_0=y=:\Phi_0(x),\label{EQ29}
\end{align}
where
\begin{align}
g^a_t:=\big(\lambda \u_{T-t}+a_t\cdot\nabla^{(2)}\Phi_t\big)\circ\Phi^{-1}_t,\ \ \Theta_t:= \big(\nn^{(2)}_{\si_t\cdot}\Phi_t)\circ\Phi_t^{-1}.\label{NB3}
\end{align}
By \eqref{NB4},  \eqref{JG1}, \eqref{ED13} and  $\si_t\in C_b^1$ uniformly in $t\in [0,T]$,   there is a constant $C>0$ such that for all $t\in[0,T]$ and $y,y'\in\mR^d$,
\begin{align}
|g^a_t(y)-g^a_t(y')|\leq C(H^{\eps'}\circ\Phi^{-1}_t(y)+H^{\eps'}\circ\Phi^{-1}_t(y'))|y-y'|,\ \ \|\nn \Theta_t\|_\infty \le C.\label{LG2}
\end{align}
On the other hand, by \eqref{1.8} and $\eps'<\eps$, for any $K>0$, there exists $C_K>0$ such that
\begin{align}
\mE\exp\bigg[K\sup_{t\in[0,T]}(H\circ\Phi^{-1}_t(Y_t))^{\eps'}\bigg]\leq C_K\exp[H(x)^{\eps}].\label{LG1}
\end{align}

In order to show the diffeomorphism property of $x\mapsto X_t(x)$, we shall use Kunita's argument. More precisely, we want to show the following estimates:
for any $p\in\mR$ and $T>0$, there are constants $C_1, C_2>0$ such that for all $x,x'\in\mR^{d_1+d_2}$ and $t\in[0,T]$,
\begin{align}
\mE|X_t(x)-X_t(x')|^{2p}&\leq C_1(\e^{H(x)^{\eps}}+\e^{H(x')^{\eps}})|x-x'|^{2p}, \label{ES1}\\
\mE(1+|X_t(x)|^{\delta_2})^{p}&\leq C_2(1+|x|^{\delta_1})^{p},\ \ p<0;\label{ES2}
\end{align}
and for any $p\geq 1$ and $T>0$, there is a constant $C_3>0$ such that for all $x\in\mR^{d_1+d_2}$ and $t,s\in[0,T]$,
\begin{align}\label{ES3}
\mE|X_t(x)-X_s(x)|^{2p}\leq C_3\e^{H(x)^{\eps}}|t-s|^{p}.
\end{align}

Estimate \eqref{ES3} is direct by the assumptions, \eqref{SDE1} and \eqref{LG1}. Let us show \eqref{ES1}. Set
$$
Z_t:=Y_t(y)-Y_t(y'),\ \ y=\Phi_0(x),\ y'=\Phi_0(x'),
$$
and
$$
G_t:=g^a_t(Y_t(y))-g^a_t(Y_t(y')),\ \ U_t:=\Theta_t(Y_t(y))-\Theta_t(Y_t(y').
$$
By It\^o's formula, we have
\begin{align*}
\dif |Z_t|^2&=[2\<Z_t,G_t\>+\tr(U^*_tU_t)]\dif t+2\<Z_t,U_t\dif W_t\>=|Z_t|^2\dif(N_t+M_t),
\end{align*}
where
$$
N_t:=\int^t_0|Z_s|^{-2}[2\<Z_s,G_s\>+\tr(U^*_sU_s)]\dif s,\ M_t:=2\int^t_0|Z_s|^{-2}\<Z_s,U_s\dif W_s\>.
$$
Here we use the convention $\frac{0}{0}=0$.
Notice that by \eqref{LG2},
\begin{align}\label{LG6}
|G_t|\leq C(H^{\eps'}\circ\Phi^{-1}_t(Y_t(y))+H^{\eps'}\circ\Phi^{-1}_t(Y_t(y')))|Z_t|,\ \ |U_t|\leq C|Z_t|.
\end{align}
Hence, by \eqref{LG1}, $N_t+M_t$ is a continuous semimartingale, and
$$
|Z_t|^2=|Z_0|^2\exp\Big\{M_t-\tfrac{1}{2}\<M\>_t+N_t\Big\}.
$$
Since for any $q\in\mR$, $t\mapsto \exp\Big\{q M_t-\tfrac{q^2}{2}\<M\>_t\Big\}$ is an exponential martingale,
by \eqref{LG6}, \eqref{LG1} and using H\"older's inequality, we have for any $p\in\mR$,
\begin{align*}
\mE |Z_t|^{2p}&=|Z_0|^{2p}\mE\exp\Big\{p M_t-\tfrac{p}{2}\<M\>_t+pN_t\Big\}\leq C(\e^{H(x)^{\eps}}+\e^{H(x')^{\eps}})|Z_0|^{2p},
\end{align*}
which in turn gives \eqref{ES1}.

Next comes to \eqref{Es2}. By It\^o's formula and \eqref{1.20}, we have
\begin{align*}
\mE H(X_t(x))^{p}=&H(x)^{p}+p\mE\int^t_0H(X_s(x))^{p-1}(\sL^{\Sigma,b+a}_s H)(X_s(x))\dif s\\
&+\frac{p(p-1)}{2}\mE\int^t_0H(X_s(x))^{p-2}|\sigma_t\cdot\nabla^{(2)}H(X_s(x))|^2\dif s\\
\leq&H(x)^{p}+C\mE\int^t_0H(X_s(x))^p\dif s,
\end{align*}
which in turn gives \eqref{ES2} by Gronwall's inequality and \eqref{1.20}.

Finally, by \eqref{ES1}-\eqref{ES3}, as in the proof of Kunita \cite[p.159-160]{Ku} (see also \cite[Theorem 3.4]{Zh0}), there is a full set $\Omega_0$
such that for all $\omega\in\Omega_0$ and $t>0$, $x\mapsto X_t(x,\omega)$ is a homeomorphism. On the other hand, since the coefficients of SDE \eqref{EQ29}
are $C^{1+\eps}$, by \cite[Theorem 4.7.2]{Ku}, $\{Y_t(\cdot)\}_{t\geq 0}$ defines a local $C^1$-diffeomorphism flow, so does
$\{X_t(\cdot)\}_{t\geq 0}$. This together with the homeomorphism property implies the global $C^1$-diffeomorphism property of $\{X_t(\cdot)\}_{t\geq 0}$. Finally,  \eqref{FL'}
  follows from \eqref{ES1} and \cite[Lemma 2.1]{Xi-Zh}.
\end{proof}

\end{document}